\documentclass[11pt]{article}

\usepackage{tikz-cd}

\usepackage[utf8]{inputenc}
\usepackage[OT2]{fontenc}
\usepackage[T1]{fontenc}

\usepackage{palatino}

\usepackage[english]{babel}
\usepackage{comment}
\usepackage{amsmath}
\usepackage{amsfonts}
\usepackage{mathrsfs}
\usepackage{amssymb}
\usepackage{geometry}
\usepackage{stmaryrd}
\usepackage{chemfig}
\usepackage{mathtools}
\usepackage{amsthm}
\usepackage{hyperref}
\usepackage{wasysym}
\usepackage{xcolor}
\usepackage{caption}
\usepackage{enumitem}
\usepackage{tikz-cd}
\usepackage{float}
\usepackage{xcolor}
\usetikzlibrary{decorations.markings}
\usetikzlibrary{shapes,positioning,intersections,quotes}
\usetikzlibrary{shapes,fit}
\usetikzlibrary{calc}

\theoremstyle{plain}
\newtheorem{theo}{Theorem}[section]
\newtheorem{prop}[theo]{Proposition}

\newtheorem{lemma}[theo]{Lemma}
\newtheorem{conj}[theo]{Conjecture}

\theoremstyle{definition}
\newtheorem{madef}[theo]{Definition}

\newtheorem*{nota*}{Notation}
\newtheorem{rmq}[theo]{Remark}
\newtheorem{ex}[theo]{Example}

\newtheorem{madeflemma}[theo]{Definition/Lemma}
\newtheorem{madefprop}[theo]{Definition/Proposition}

\renewcommand{\L}{\mathcal{L}}
\newcommand{\zero}{\hat{0}}
\newcommand{\rk}{\mathrm{rk}\,}
\newcommand{\medvee}{\scalebox{1.1}{$\bigvee$}}
\newcommand{\Z}{\mathbb{Z}}
\newcommand{\E}{\mathcal{E}}
\newcommand{\MD}{\mathcal{MD}}
\newcommand{\isomod}{\vee}
\newcommand{\un}{\hat{1}}
\renewcommand{\d}{\mathrm{d}}
\newcommand{\At}{\mathrm{At}}
\newcommand{\OS}{\mathrm{OS}}
\newcommand{\I}{\mathbb{I}}
\newcommand{\Ibar}{\overline{\I}}
\newcommand{\Sc}{J^{\mathrm{contr}}}
\newcommand{\Id}{\mathrm{Id}}
\newcommand{\Q}{\mathbb{Q}}
\newcommand{\Cx}{\mathbb{C}}
\newcommand{\A}{\mathcal{A}}
\newcommand{\Br}{\mathcal{B}r}
\newcommand{\M}{\mathcal{M}}
\newcommand{\Tr}{\mathrm{Tr}}
\newcommand{\Sym}{\mathbb{S}}
\newcommand{\gr}{\mathrm{gr}\,}
\newcommand{\Qu}{\mathcal{Q}}
\renewcommand{\H}{\mathrm{H}}

\renewcommand{\C}{\mathcal{C}}
\newcommand{\GL}{\mathcal{GL}}
\newcommand{\SFM}{\mathrm{SFM}}
\newcommand{\R}{\mathbb{R}}
\newcommand{\Gra}{\textbf{Gra}}
\renewcommand{\S}{\mathrm{S}}
\newcommand{\Conf}{\mathrm{Conf}}
\newcommand{\Harr}{\mathcal{H}}
\newcommand{\qA}{\widetilde{\A}}
\renewcommand{\P}{\mathbb{P}}
\newcommand{\FM}{\mathrm{FM}}
\newcommand{\G}{\mathcal{G}}
\newcommand{\codim}{\mathrm{codim}}
\newcommand{\vol}{\mathrm{vol}}
\newcommand{\pr}{\mathrm{pr}}
\newcommand{\f}{\widetilde{f}}
\newcommand{\pibar}{\overline{\pi}}
\newcommand{\bigarr}{\mathcal{K}}
\newcommand{\h}{\mathfrak{h}}
\newcommand{\Lie}{\mathfrak{Lie}}
\newcommand{\IMD}{\mathcal{IMD}}
\newcommand{\F}{\mathfrak{F}}
\newcommand{\LD}{\mathrm{LD}}
\newcommand{\N}{\mathbb{N}}
\newcommand{\opi}{\overline{\pi}}

\title{Matroid complexes and Orlik--Solomon algebras}
\author{Basile Coron}
\date{}

\begin{document}
\maketitle

\begin{abstract}
In this article we construct a combinatorial quasi-free differential graded model for the Orlik--Solomon algebra of supersolvable matroids, which generalizes in a matroidal setting the cdga of admissible graphs introduced by M. Kontsevich for the braid arrangements. Our construction draws on well-known concepts from matroid theory, including modularity, single-element extensions, and generalized parallel connections. We also show that this model carries a cooperadic structure in a suitably generalized sense. As an application, we use this model to give a new proof that the Orlik–Solomon algebras of supersolvable matroids are Koszul. 
\end{abstract}

\section{Introduction}
The little $2$-disk operad $\LD$ is the object consisting of the collection $\{\LD(n) \, | \, n \in \mathbb{N}\}$ of spaces of configurations of $n$ disjoint $2$-disks inside the unit $2$-disk, together with the morphisms $$\LD(p)\times \LD(q) \rightarrow \LD(p+q-1)\, \textrm{ for all } i\leq p,q$$ given by inserting the second configuration in the $i$-th disk of the first one. It was originally introduced by Boardman and Vogt \cite{BV_1968} and May \cite{May_1972} to study iterated loop spaces. Since then, it has found applications in other areas, such as Goodwillie–Weiss embedding calculus \cite{GW_1999, Weiss_1999} and deformation-quantization of Poisson manifolds \cite{Kontsevich1999, Kontsevich2003} in mathematical physics. For a general overview we refer to Fresse \cite[Section 11]{Miller_2020}. Each space $\LD(n)$ is homotopy equivalent to the configuration space $\Conf_n(\Cx)$ of $n$ distinct points in the complex plane. However, unlike $\LD(n)$, these configuration spaces do not support a natural insertion operation. One alternative to obtaining such a structure—besides “fattening the points”—is to consider the (spherical) Fulton--MacPherson compactifications $\FM(n) \supset \Conf_{n}(\Cx)$ as introduced by Axelrod and Singer \cite{Axelrod_Singer_1993} or Fulton and MacPherson \cite{Fulton1994}. The inclusion of strata of the spaces $\FM(n)$ give morphisms $$\FM(p)\times \FM(q) \rightarrow \FM(p+q-1)\, \textrm{ for all }  i\leq p,q$$ which constitute the Fulton--MacPherson operad $\FM$. By a result of Salvatore \cite{Salvatore_2001} the operads $\LD$ and $\FM$ are equivalent. 

The rational cohomology rings $\H^{\bullet}(\Conf_n(\Cx), \Q) \simeq \H^{\bullet}(\FM(n), \Q)$ admit the following explicit presentation due to Arnold \cite{Arnold1969}: 
\begin{equation*}
    \H^{\bullet}(\Conf_n(\Cx), \Q) \simeq \frac{\bigwedge^{\bullet}[e_{ij} \, | \,1 \leq i<j\leq n , \deg(e_{ij})=1]}{( e_{ij}e_{ik} - e_{ij}e_{jk} + e_{ik}e_{jk}, i<j<k)},
\end{equation*}
where each $e_{ij}$ corresponds to the differential $1$-form $\omega_{ij} = \frac{d(z_i - z_j)}{z_i - z_j}$. To study the rational homotopy type of $\Conf_n(\Cx)$, rational homotopy theory \cite{Quillen_1969, Sullivan1977} teaches us to look for a quasi-free model of $\H^{\bullet}(\Conf_n(\Cx), \Q)$, that is, a commutative differential graded algebra (cdga for short) which is free as an algebra and which is quasi-isomorphic to $\H^{\bullet}(\Conf_n(\Cx), \Q)$. In principle, one can always construct such a model inductively: starting with the free graded commutative algebra $\bigwedge^{\bullet}[e_{ij} \, | \,i<j, \deg(e_{ij})=1]$, one can add formal generators $r_{ijk}$ with differential $\d(r_{ijk})= e_{ij}e_{ik} - e_{ij}e_{jk} + e_{ik}e_{jk}$, and then continue adding generators to resolve all remaining relations. However, in practice, constructing a manageable and explicit quasi-free model can be quite challenging. In \cite{Kontsevich2003} Kontsevich introduced such a model: the cdga of admissible graphs $\Gra(n)$ which is quasi-isomorphic to $\H^{\bullet}(\FM(n), \Q)$. That cdga has a linear basis given by certain graphs with $n$ numbered "external" vertices and other additional unnumbered "internal" ones. The differential is defined as an alternating sum of edge contractions, and the product is given by gluing graphs along their external vertices. Moreover, the collection $\{\Gra(n)\, |\, n\in \N\}$ carries a (co)operadic structure given by morphisms of cdgas
$$\Gra(p+q-1) \rightarrow \Gra(p)\otimes \Gra(q) \, \textrm{ for all } i\leq p, q$$ 
which reflect the operadic structure of $\FM$. Each graph $\Gamma \in \Gra(n)$ with $k$ external vertices corresponds to an explicit piecewise algebraic form $$\int_{\FM(n+k)\twoheadrightarrow \FM(n)}\bigwedge_{(i,j) \textrm{ edge of } \Gamma}\omega_{i,j} \in \Omega^{\bullet}_{PA}(\FM(n), \R),$$ where $\FM(n+k)\twoheadrightarrow \FM(n)$ is the fibration given by forgetting the last $k$ points, and the integral denotes the pushforward along this fibration. This construction yields a zig-zag of qisos of cooperads of cdgas (modulo some technicalities)
\begin{equation*}
    \Omega_{PA}^{\bullet}(\FM(n), \R) \xleftarrow{\sim} \Gra(n)\otimes \R \xrightarrow{\sim} \H^{\bullet}(\FM(n), \R). 
\end{equation*}
The explicit construction of $\Gra(n)$ is briefly recalled in Section \ref{secgra}. We refer to \cite{Kontsevich2003} and \cite{lambrechts_formality_2012} for more details. More recently, Idrissi \cite{Idrissi2019} extended the graphical model $\Gra$ to configuration spaces of manifolds beyond $\Cx$. In this paper, we explore a different direction of generalization. The configuration space $\Conf_n(\Cx)$ can be viewed as the complement of the braid arrangement, which is the collection of hyperplanes $\{z_i = z_j\}$ in $\Cx^n$. Our main goal is to investigate to what extent the above story can be extended to arbitrary complexified hyperplane arrangements, and even further, to arbitrary matroids, which abstract the combinatorics of such arrangements (see \cite{OT_1992} for a general reference on hyperplane arrangements and \cite{oxley} for a general reference on matroids).

In Section \ref{secMD}, we associate a cdga $\MD(\L)$ to any matroid $\L$, which has a basis given by pairs $(\E, J)$ of a matroid $\E$ containing $\L$ as a modular flat and a subset $J$ of atoms of $\E$. The differential of an element $(\E, J)$ is given by the alternating sum of the contractions of $\E$ along the elements of $J$, and the commutative product is given by taking the generalized parallel connection along $\L$. The geometric motivations behind the definition of $\MD(\L)$ are discussed in Section \ref{secgeo}. In short if $\L$ is realized by some complexified hyperplane arrangement $\Harr$, an element $(\E, J) \in \MD(\L)$ with $\E$ realized by some larger complexified hyperplane arrangement $\bigarr$ corresponds to a particular piecewise algebraic form on the (spherical) Fulton--Macpherson compactification $\FM[\Harr]$ of the complement of $\Harr$, as introduced by Gaiffi \cite{Gaiffi_2003}. That form is defined as the pushforward 
$$\int_{\FM[\bigarr] \twoheadrightarrow \FM[\Harr]}\bigwedge_{H \in J} \d f_H/f_H \in \Omega_{PA}^{\bullet}(\FM[\Harr], \R)$$ along the restriction fibration $\FM[\bigarr] \twoheadrightarrow \FM[\Harr]$ (where the $f_H$'s are some chosen annihilators of the hyperplanes of $\bigarr$). For an arbitrary matroid $\L$, Proposition \ref{propmor} gives an explicit morphism  of cdgas $\I_{\L}:\MD(\L) \rightarrow \OS(\L)$, where $\OS(\L)$ denotes the Orlik--Solomon algebra of $\L$, which is the cohomology algebra of the arrangement complement if $\L$ is realizable over $\Cx$ (see Orlik and Solomon \cite{OS_1980}). In general $\I_{\L}$ is not a quasi-isomorphism because $\MD(\L)$ only resolves the quadratic relations of $\OS(\L)$. The braid arrangements have the particular property that we have a sequence of fibrations between the associated arrangement complements 
\begin{equation*}
    \Conf_n(\Cx) \twoheadrightarrow \Conf_{n-1}(\Cx) \twoheadrightarrow \cdots \twoheadrightarrow \Cx^{\star} \twoheadrightarrow \{\star\},
\end{equation*}
given by forgetting points one after the other. In general an arrangement having such a sequence of fibrations is called fiber-type. For instance Coxeter arrangements in type B are also fiber-type. Combinatorially, the fiber-type condition corresponds to the supersolvability of the underlying matroid, a notion introduced by Stanley \cite{Stanley_1972}. Theorem \ref{theoqiso} shows that if a matroid $\L$ is supersolvable then the morphism $\I_{\L}: \MD(\L) \rightarrow \OS(\L)$ is a quasi-isomorphism. 

In Section \ref{secstruct} we construct morphisms of cdgas $$\MD(\L) \rightarrow \MD(\L_F)\otimes \MD(\L^F), \textrm{ for all matroid } \L \textrm{ and all } F \in \L,$$
where $\L_F$ and $\L^F$ denote the restriction and contraction of $\L$ at $F$ respectively. This defines an object called a $\GL$-(co)operad, as defined by the author in \cite{BC2025,BC2024}. Those morphisms model the inclusions of strata 
$$\FM[\Harr_F]\times \FM[\Harr^F] \rightarrow \FM[\Harr]$$
in the realizable case. If we restrict to braid arrangements this is the classical operadic structure on $\FM$ mentioned above.  

Finally, in Section \ref{seckos} we show how the cdga $\MD(\L)$ can be used to give a new proof of the Koszulness of Orlik--Solomon algebras of supersolvable matroids, generalizing a strategy suggested in \cite{SW_2009} for the braid case. \\

\textbf{Acknowledgment.} This research was carried out while the author was at the Institute for Advanced Study in Princeton, on the occasion of the special year program in combinatorial geometry. The author is grateful to the organizers and the participants of the program for an enjoyable semester. 
\section{The cdga of admissible graphs}\label{secgra}
The cdga of admissible graphs was introduced by Kontsevich \cite{Kontsevich2003} as a combinatorial dg real model of configuration spaces of $\R^2$. The construction was later worked out in details by Lambrechts and Volic in \cite{lambrechts_formality_2012}. In this section we recall the main definitions of that story, following the exposition of \cite{lambrechts_formality_2012} and referring to that article for more details. Note that there is also an odd version of this cdga (modeling configuration spaces of $\R^3$ instead of $\R^2$) but in this article we will exclusively consider the even case. The natural home for the odd version would be oriented matroids instead of matroids. \\

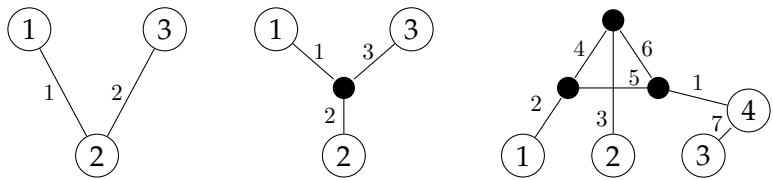
\begin{figure}[H]
    \centering
    \begin{tikzpicture}[scale=0.60, baseline=14px]
    \node[] (A) at (0.5,2.8) {1};
    \node[] (B) at (2,0) {2};
    \node[] (C) at (3.5,2.8) {3};
    \node[circle,draw, inner xsep=-0.5mm,inner ysep=-0.5mm, fit= (A)] (ca) {};
    \node[circle,draw, inner xsep=-0.5mm,inner ysep=-0.5mm, fit= (B)] (cb) {};
    \node[circle,draw, inner xsep=-0.5mm,inner ysep=-0.5mm, fit= (C)] (cc) {};
    
    \draw[] (ca) -- (cb) node[midway,left, scale = 0.8] {$1$};
    \draw[] (cc) -- (cb) node[midway,left, scale = 0.8] {$2$};
    \end{tikzpicture}
    \hspace{20pt}
\begin{tikzpicture}[scale=0.60, baseline=14px]
    
    \fill[black] (2,1.5) circle (7pt);
    \node[] (A) at (0.5,2.8) {1};
    \node[] (B) at (2,0) {2};
    \node[] (C) at (3.5,2.8) {3};
    \node[] (D) at (2,1.5) {}; 
    \node[circle,draw, inner xsep=-0.5mm,inner ysep=-0.5mm, fit= (A)] (ca) {};
    \node[circle,draw, inner xsep=-0.5mm,inner ysep=-0.5mm, fit= (B)] (cb) {};
    \node[circle,draw, inner xsep=-0.5mm,inner ysep=-0.5mm, fit= (C)] (cc) {};
    
    \draw[] (ca) -- (D) node[pos=0.65,above, scale = 0.8] {$1$};
    \draw[] (cb) -- (D) node[midway,left, scale = 0.8] {$2$};
    \draw[] (cc) -- (D) node[pos=0.65,above, scale = 0.8] {$3$};
\end{tikzpicture}
\hspace{20pt}
\begin{tikzpicture}[scale=0.60, baseline=14px]
    \fill[black] (2,3) circle (7pt);
    \fill[black] (1,1.5) circle (7pt);
    \fill[black] (3,1.5) circle (7pt);
    \node[] (A) at (0,0) {1};
    \node[] (B) at (2,0) {2};
    \node[] (C) at (4,0) {3};
    \node[] (D) at (5,1) {4};
    \node[] (E) at (2,3) {}; 
    \node[] (F) at (1,1.5) {}; 
    \node[] (G) at (3,1.5) {}; 
    \node[circle,draw, inner xsep=-0.5mm,inner ysep=-0.5mm, fit= (A)] (ca) {};
    \node[circle,draw, inner xsep=-0.5mm,inner ysep=-0.5mm, fit= (B)] (cb) {};
    \node[circle,draw, inner xsep=-0.5mm,inner ysep=-0.5mm, fit= (C)] (cc) {};
    \node[circle,draw, inner xsep=-0.5mm,inner ysep=-0.5mm, fit= (D)] (cd) {};
    \draw[] (cd) -- (G) node[midway,above, scale = 0.8] {$1$};
    \draw[] (ca) -- (F) node[midway,above left=-1pt, scale = 0.8] {$2$};
    \draw[] (cb) -- (E) node[pos=0.15 , left=-1pt, scale = 0.8] {$3$};
    \draw[] (E) -- (F)  node[pos=0.6, above left = -3pt, scale = 0.8] {$4$};
    \draw[] (F) -- (G)  node[pos=0.8,above = -2pt, scale = 0.8] {$5$};
    \draw[] (G) -- (E)  node[pos=0.4,above right = -3pt, scale = 0.8] {$6$} ;
    \draw[] (cc) -- (D) node[midway,above left = -3pt, scale = 0.8] {$7$};
    
\end{tikzpicture}
    \caption{Examples of diagrams}
    \label{exdiag}
\end{figure}

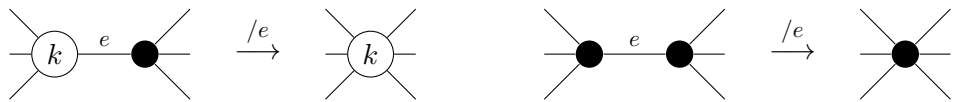
\begin{figure}[H]
    \centering
    \begin{tikzpicture}[scale=0.60, baseline=14px]
   
    \node[circle,draw, inner xsep=1mm,inner ysep=0.9mm] (ca) at (0,2) {$k$};
    \node[circle,fill] (cb) at (2,2) {};
    
    \draw[] (ca) -- (cb) node[midway,above, scale = 0.8] {$e$};
    \draw[] (ca) -- (-1,3); 
    \draw[] (ca) -- (-1,2); 
    \draw[] (ca) -- (-1,1); 
    \draw[] (cb) -- (3,3); 
    \draw[] (cb) -- (3,2); 
    \draw[] (cb) -- (3,1); 
    \node[] (arr) at (4.5,2) {$\longrightarrow$};
    \node[scale=0.8] (arr) at (4.5,2.5) {$/e$};
    \node[circle,draw, inner xsep=1mm,inner ysep=0.9mm] (cc) at (7,2) {$k$};
    \draw[] (cc) -- (6,3); 
    \draw[] (cc) -- (6,2); 
    \draw[] (cc) -- (6,1); 
    \draw[] (cc) -- (8,3); 
    \draw[] (cc) -- (8,2); 
    \draw[] (cc) -- (8,1);

    \end{tikzpicture}
    \hspace{20pt}
\hspace{20pt}
\begin{tikzpicture}[scale=0.60, baseline=14px]
   
    \node[circle,fill] (ca) at (0,2) {};
    \node[circle,fill] (cb) at (2,2) {};
    \draw[] (ca) -- (-1,3); 
    \draw[] (ca) -- (-1,2); 
    \draw[] (ca) -- (-1,1); 
    \draw[] (cb) -- (3,3); 
    \draw[] (cb) -- (3,2); 
    \draw[] (cb) -- (3,1); 
    
    \draw[] (ca) -- (cb) node[midway,above, scale = 0.8] {$e$};

    \node[] (arr) at (4.5,2) {$\longrightarrow$};
    \node[scale=0.8] (arr) at (4.5,2.5) {$/e$};
    \node[circle,fill] (cc) at (7,2) {};
    \draw[] (cc) -- (6,3); 
    \draw[] (cc) -- (6,2); 
    \draw[] (cc) -- (6,1); 
    \draw[] (cc) -- (8,3); 
    \draw[] (cc) -- (8,2); 
    \draw[] (cc) -- (8,1);
   
    \end{tikzpicture}
    \caption{Contracting a contractible edge}
    \label{contraction}
\end{figure}

For any integer $n\geq 1$, a \textit{diagram of arity $n$} is a non-oriented loopless graph with vertices of two kinds: the external vertices, which are numbered from $1$ to $n$, and the internal vertices which are unnumbered and can have any cardinality. A diagram also comes equipped with a total order on its edges. See Figure \ref{exdiag} for some examples. 
A diagram is called \textit{admissible} if every internal vertex has at least $3$ neighbours, and every internal vertex has an external vertex in its connected component. We denote by $\Gra(n)$ the $\Q$-vector space spanned by diagrams of arity $n$ quotiented by the relations $\Gamma \sim 0$ if $\Gamma$ is not admissible, and $\Gamma\sim -\Gamma'$ if $\Gamma$ and $\Gamma'$ differ by a transposition of the linear order on their edges. The vector space $\Gra(n)$ is $\Z$-graded by setting the degree of a diagram $\Gamma$ to be $-2\#inner(\Gamma) + \#edges(\Gamma)$ where $\#inner(\Gamma)$ denotes the number of internal vertices of $\Gamma$ and $\#edges(\Gamma)$ denotes its number of edges. An edge of $\Gamma$ is called contractible if it is not an edge between two external vertices and if no extremity of $e$ is a univalent internal vertex. For any contractible edge $e$ of a diagram $\Gamma$, one can contract $e$ to get another diagram $\Gamma/e$, as described locally in Figure \ref{contraction}. 
\vspace{10pt}
The graded vector space $\Gra(n)$ admits a differential $\d$ defined as the alternating sum of all possible contractions: 
\begin{equation*}
    \d(\Gamma) \coloneqq \sum_{e \textrm{ contractible}}\pm \Gamma/e.
\end{equation*}
For instance we have 
\begin{equation}\label{exdiff}
\d
    \begin{tikzpicture}[scale=0.60, baseline=14px]
    
    \fill[black] (2,1) circle (7pt);
    \node[] (A) at (0.5,2.3) {1};
    \node[] (B) at (2,-0.5) {2};
    \node[] (C) at (3.5,2.3) {3};
    \node[] (D) at (2,1) {}; 
    \node[circle,draw, inner xsep=-0.5mm,inner ysep=-0.5mm, fit= (A)] (ca) {};
    \node[circle,draw, inner xsep=-0.5mm,inner ysep=-0.5mm, fit= (B)] (cb) {};
    \node[circle,draw, inner xsep=-0.5mm,inner ysep=-0.5mm, fit= (C)] (cc) {};
    
    \draw[] (ca) -- (D) node[pos=0.65,above, scale = 0.8] {$1$};
    \draw[] (cb) -- (D) node[midway,left, scale = 0.8] {$2$};
    \draw[] (cc) -- (D) node[pos=0.65,above, scale = 0.8] {$3$};
\end{tikzpicture} 
= 
\begin{tikzpicture}[scale=0.60, baseline=14px]
    \node[] (A) at (0.5,2.3) {1};
    \node[] (B) at (2,-0.5) {2};
    \node[] (C) at (3.5,2.3) {3};
    \node[circle,draw, inner xsep=-0.5mm,inner ysep=-0.5mm, fit= (A)] (ca) {};
    \node[circle,draw, inner xsep=-0.5mm,inner ysep=-0.5mm, fit= (B)] (cb) {};
    \node[circle,draw, inner xsep=-0.5mm,inner ysep=-0.5mm, fit= (C)] (cc) {};
    
    \draw[] (ca) -- (cc) node[midway,above, scale = 0.8] {$2$};
    \draw[] (ca) -- (cb) node[midway,left, scale = 0.8] {$1$};
    \end{tikzpicture}
    -
    \begin{tikzpicture}[scale=0.60, baseline=14px]
    \node[] (A) at (0.5,2.3) {1};
    \node[] (B) at (2,-0.5) {2};
    \node[] (C) at (3.5,2.3) {3};
    \node[circle,draw, inner xsep=-0.5mm,inner ysep=-0.5mm, fit= (A)] (ca) {};
    \node[circle,draw, inner xsep=-0.5mm,inner ysep=-0.5mm, fit= (B)] (cb) {};
    \node[circle,draw, inner xsep=-0.5mm,inner ysep=-0.5mm, fit= (C)] (cc) {};
    
    \draw[] (ca) -- (cb) node[midway,left, scale = 0.8] {$1$};
    \draw[] (cc) -- (cb) node[midway,left, scale = 0.8] {$2$};
    \end{tikzpicture}
    +
    \begin{tikzpicture}[scale=0.60, baseline=14px]
    \node[] (A) at (0.5,2.3) {1};
    \node[] (B) at (2,-0.5) {2};
    \node[] (C) at (3.5,2.3) {3};
    \node[circle,draw, inner xsep=-0.5mm,inner ysep=-0.5mm, fit= (A)] (ca) {};
    \node[circle,draw, inner xsep=-0.5mm,inner ysep=-0.5mm, fit= (B)] (cb) {};
    \node[circle,draw, inner xsep=-0.5mm,inner ysep=-0.5mm, fit= (C)] (cc) {};
    
    \draw[] (ca) -- (cc) node[midway,above, scale = 0.8] {$1$};
    \draw[] (cc) -- (cb) node[midway,left, scale = 0.8] {$2$};
    \end{tikzpicture}
\end{equation}
The graded vector space $\Gra(n)$ admits a graded commutative product $\bullet$ defined on admissible diagrams by gluing graphs along the external vertices and concatenating the orders on the edges. For instance we have 
\begin{equation*}
\begin{tikzpicture}[scale=0.60, baseline=14px]
    
    \fill[black] (2,1) circle (7pt);
    \node[] (A) at (0.5,2.3) {1};
    \node[] (B) at (2,-0.5) {2};
    \node[] (C) at (3.5,2.3) {3};
    \node[] (D) at (2,1) {}; 
    \node[circle,draw, inner xsep=-0.5mm,inner ysep=-0.5mm, fit= (A)] (ca) {};
    \node[circle,draw, inner xsep=-0.5mm,inner ysep=-0.5mm, fit= (B)] (cb) {};
    \node[circle,draw, inner xsep=-0.5mm,inner ysep=-0.5mm, fit= (C)] (cc) {};
    
    \draw[] (ca) -- (D) node[pos=0.65,above, scale = 0.8] {$1$};
    \draw[] (cb) -- (D) node[midway,left, scale = 0.8] {$2$};
    \draw[] (cc) -- (D) node[pos=0.65,above, scale = 0.8] {$3$};
    \draw[] (ca) -- (cc) node[midway,above, scale = 0.8] {$4$};
\end{tikzpicture} 
\bullet 
\begin{tikzpicture}[scale=0.60, baseline=14px]
    \node[] (A) at (0.5,2.3) {1};
    \node[] (B) at (2,-0.5) {2};
    \node[] (C) at (3.5,2.3) {3};
    \node[circle,draw, inner xsep=-0.5mm,inner ysep=-0.5mm, fit= (A)] (ca) {};
    \node[circle,draw, inner xsep=-0.5mm,inner ysep=-0.5mm, fit= (B)] (cb) {};
    \node[circle,draw, inner xsep=-0.5mm,inner ysep=-0.5mm, fit= (C)] (cc) {};
    
    \draw[] (ca) -- (cb) node[midway,left, scale = 0.8] {$1$};
    \draw[] (cc) -- (cb) node[midway,left, scale = 0.8] {$2$};
    \end{tikzpicture}
=
\begin{tikzpicture}[scale=0.60, baseline=14px]
    
    \fill[black] (2,1) circle (7pt);
    \node[] (A) at (0.5,2.3) {1};
    \node[] (B) at (2,-0.5) {2};
    \node[] (C) at (3.5,2.3) {3};
    \node[] (D) at (2,1) {}; 
    \node[circle,draw, inner xsep=-0.5mm,inner ysep=-0.5mm, fit= (A)] (ca) {};
    \node[circle,draw, inner xsep=-0.5mm,inner ysep=-0.5mm, fit= (B)] (cb) {};
    \node[circle,draw, inner xsep=-0.5mm,inner ysep=-0.5mm, fit= (C)] (cc) {};
    
    \draw[] (ca) -- (D) node[pos=0.65,above, scale = 0.8] {$1$};
    \draw[] (cb) -- (D) node[midway,left = -2pt, scale = 0.8] {$2$};
    \draw[] (cc) -- (D) node[pos=0.65,above, scale = 0.8] {$3$};
    \draw[] (ca) -- (cc) node[midway,above, scale = 0.8] {$4$};
    \draw[] (ca) -- (cb) node[pos=0.3,below = 2pt, scale = 0.8] {$5$};
    \draw[] (cb) -- (cc) node[pos=0.7,below = 2pt, scale = 0.8] {$6$};
\end{tikzpicture} 
\end{equation*}
One can see that the commutative algebra $(\Gra(n), \bullet)$ is free with generators given by admissible diagrams which are internally connected, that is which are connected after removing the external vertices. The commutative product $\bullet$ is compatible with the differential, in the sense that it satisfies the Leibniz identity
\begin{equation*}
    \d(\Gamma_1\bullet \Gamma_2) = \d(\Gamma_1)\bullet \Gamma_2 + (-1)^{\deg(\Gamma_1)}\Gamma_1\bullet \d(\Gamma_2),
\end{equation*}
which makes $(\Gra(n), \bullet, \d)$ a cdga. Kontsevich \cite{Kontsevich1999} proved that the cdga $\Gra(n)$ is a model of the rational cohomology algebra of $\Conf_n(\Cx)$, meaning that we have a quasi-isomorphism of cdgas $\Gra(n)\xrightarrow{\sim}\H^{\bullet}(\Conf_n(\Cx), \Q)$. To describe this quasi-isomorphism let us recall that we have an explicit presentation of $\H^{\bullet}(\Conf_n(\Cx), \Q)$ due to Arnold \cite{Arnold1969}: 
\begin{equation}\label{presArnold}
    \H^{\bullet}(\Conf_n(\Cx), \Q) \simeq \frac{\bigwedge^{\bullet}[e_{ij} \, | \, 1 \leq i<j\leq n, \deg(e_{ij})=1]}{( e_{ij}e_{ik} - e_{ij}e_{jk} + e_{ik}e_{jk}, i<j<k)},
\end{equation}
where $\bigwedge^{\bullet}$ denotes the free graded commutative algebra. The following explicit quasi-isomorphism is due to Lambrechts and Volic. 
\begin{theo}[{\cite[Theorem 8.1]{lambrechts_formality_2012}}]
There is a quasi-isomorphism of cdgas $\Ibar: \Gra(n) \xrightarrow{\sim} (\H^{\bullet}(\Conf_n(\Cx),\Q), \d=0)$ such that $\Ibar$ is zero on admissible diagrams containing a contractible edge, and such that we have 
\begin{equation*}
    \Ibar( 
\begin{tikzpicture}[scale=0.60, baseline=14px]
    \node[] (A) at (0.5,0.97) {i};
    \node[scale=0.88] (B) at (1.8,0.97) {j};
    \node[circle,draw, inner xsep=-3pt,inner ysep=-3pt, fit= (A)] (ca) {};
    \node[circle,draw, inner xsep=-3pt,inner ysep=-3pt, fit= (B)] (cb) {};
    \draw[] (ca) -- (cb);
\end{tikzpicture}
) = e_{ij}
\end{equation*}
for all $i<j.$
\end{theo}
Note that every admissible diagram without contractible edges can be expressed as a product of diagrams with only one edge, and so the morphism of cdgas $\Ibar$ is completely determined. The fact that this morphism is well-defined means that for all admissible diagram $\Gamma$ we have $\Ibar(\d(\Gamma))=0$. For instance if $\Gamma$ is the trident in equation \eqref{exdiff}, that equation gives $\Ibar(\d(\Gamma)) = e_{12}e_{13} - e_{12}e_{23} + e_{13}e_{23}$ which is indeed $0$ in $\H^{\bullet}(\Conf_n(\Cx), \Q)$. In other words, the trident resolves the relation $e_{12}e_{13} - e_{12}e_{23} + e_{13}e_{23}=0 $. \\

For all $n\geq 1$, the cdga $\Gra(n)$ admits an $\Sym_n$-action given by permuting the $n$ external vertices. In addition, the collection of cdgas $\{\Gra(n), \, n \in \N\}$ admits a cooperadic structure given by morphisms 
\begin{equation*}
    \Gra(n) \rightarrow \Gra(p)\otimes \Gra(q), \textrm{ for all } p,q \textrm{ s.t. } p+q = n+1 \textrm{ and } i\leq p.
\end{equation*}
Those morphisms are defined by sending a diagram $\Gamma$ to the sum $\sum_{\gamma}\Gamma/\gamma\otimes \gamma$, where $\gamma$ runs over all subgraphs of $\Gamma$ containing the external vertices $i,\ldots,i+q-1$ and no other external vertices, and $\Gamma/\gamma$ denotes the diagram obtained from $\Gamma$ by contracting $\gamma$.\\

The cdga $\Gra(n)\otimes \R$ is also quasi-isomorphic to the cdga of real semi-algebraic forms of the spherical Fulton-MacPherson compactification $\FM(n)\supset \Conf_n(\Cx)$ (see \cite[Sections 4,5]{lambrechts_formality_2012}), which is homotopy equivalent to $\Conf_n(\Cx)$. The qiso $\I: \Gra(n)\otimes \R \xrightarrow{\sim} \Omega_{PA}(\FM(n), \R)$ is defined by sending admissible diagrams to some explicit semi-algebraic forms of $\FM(n)$. The properties of those forms are what motivates the definition of $\Gra(n)$ and in particular the notion of admissibility. In the next section we will explain in more details what this means, in the more general setting of hyperplane arrangements. Note that tensoring with $\R$ is necessary because the definition of the forms $\I(\Gamma)$ involves transcendental numbers. It is remarkable however, that those transcendental numbers are not necessary to define $\Gra$.

\section{The geometry}\label{secgeo}
A \textit{hyperplane arrangement} is a finite collection of hyperplanes inside a shared finite dimensional vector space. In this paper we will focus on complexified hyperplane arrangements, that is, arrangements of the form $\{H_1 \otimes\Cx,\ldots, H_n\otimes \Cx\}$ for some real arrangement $\{H_1, \ldots, H_n\}$. For instance we have the so-called braid arrangements $\Br_n = \{\{z_i = z_j\}, 1 \leq i\neq j \leq n \}$ inside $\Cx^n$. More generally to a simple graph $\Gamma$ with vertex set $V$ and edge set $E$, one can associate the graphical hyperplane arrangement $\{\{z_i = z_j\}, (i,j) \in E\}$ in $\Cx^V$. The \textit{arrangement complement} $\A_{\Harr}$ of a hyperplane arrangement $\Harr = \{H_1, \ldots, H_k\}$ in a vector space $V$ is the algebraic variety $\A_{\Harr} \coloneqq V\setminus \bigcup_{i \leq k} H_i$. For instance for the braid arrangements we have $\A_{\Br_n} = \Conf_n(\Cx).$ A hyperplane arrangement $\{H_1, \ldots, H_k\}$ is called \textit{essential} if we have $\bigcap_{i\leq k}H_i=\{0\}$. Any hyperplane arrangement in some vector space $V$ can be essentialized by considering the induced arrangement in the quotient $V/\bigcap_{i\leq k }H_i$. For instance the braid arrangement $\Br_n$ is not essential because the intersection of the diagonal hyperplanes is the full diagonal $\{z_1 = \cdots = z_n\}$. The complements of the essentialization of the braid arrangements are the configuration spaces of $\Cx$ quotiented by translations. In the rest of this paper every hyperplane arrangement will be assumed to be essential. In particular we will now use $\Br_n$ to denote the essentialization of the braid arrangements. For any complexified hyperplane arrangement $\Harr$ let us denote $\qA_{\Harr}\coloneqq \A_{\Harr}/\R_{>0}$. In \cite{Gaiffi_2003} Gaiffi showed how to construct a semi-algebraic compactification $\SFM[\Harr]$ of $\qA_{\Harr}$, generalizing a construction of Axelrod and Singer \cite{Axelrod_Singer_1993} in the braid case. Let us recall the main details of that story, referring to \cite{Gaiffi_2003} for more details. The construction of $\SFM[\Harr]$ is analogous to that of the wonderful compactification of the projectivization $\P\A_{\Harr}$, due to De Concini and Procesi \cite{DeConcini1995}. One simply needs to replace projective spaces by spheres. The letter $\S$ stands for spherical and $\FM$ stands for Fulton--Macpherson, in reference to the classical Fulton--Macpherson compactifications of configuration spaces \cite{Fulton1994}. \\ 

For any complex vector space $V$ we denote $\S(V)\coloneqq V/\R_{>0}$ the ``spherification'' of $V$. For any hyperplane arrangement $\Harr = \{H_1,\ldots, H_k\}$ the \textit{intersection lattice} $\L_{\Harr}$ of $\Harr$ is the set $\{\bigcap_{i \in I} H_i, I\subset \{1,\ldots,k\}\}$ of all possible intersections of hyperplanes of $\Harr$ ordered by reverse inclusion. The elements of $\L_{\Harr}$ will be called \textit{flats} of $\Harr$. For instance for the braid arrangement $\Br_n$ the poset $\L_{\Br_n}$ can be identified with the collection $\Pi_n$ of set partitions of $\{1,\ldots,n\}$ ordered by refinement, also called the $n$-th partition lattice. A flat is called \textit{proper} if it is neither $V$ nor $\{0\}$. If $\Harr$ is a complexified arrangement, for any flat $F \in \L_{\Harr}$ different from $V$ we have the projection $V\twoheadrightarrow V/F$ which induces a map $\qA_{\Harr}\rightarrow \S(V/F)$. Putting those maps together we obtain an embedding $\qA_{\Harr}\hookrightarrow \prod_{F\in \L_{\Harr}\setminus \{V\}}\S(V/F)$. Let $\SFM[\Harr]$ be the closure of the image of this embedding. Gaiffi \cite{Gaiffi_2003} proved that $\SFM[\Harr]$ is a $\C^{\infty}$ real manifold with corners, that is, a space locally modeled by open subsets of positive orthants $\R_{\geq 0}^{p}$ for some $p$. In fact $\SFM[\Harr]$ is a semi-algebraic manifold in the sense of \cite{bochnak_coste_roy_1998}. As a topological space $\SFM[\Harr]$ is a homotopy retract of $\qA_{\Harr}$ and so those spaces have the same homotopy type. As in the theory of wonderful compactifications, there are variations of $\SFM[\Harr]$ where one considers more general product embeddings $\qA_{\Harr}\hookrightarrow \prod_{F \in \G}\S(V/F)$ over certain subsets $\G\subset \L_{\Harr}$ called ``building sets''. For $\SFM$ however, changing $\G$ does not affect the homotopy type, which is always that of the arrangement complement. For simplicity in this article we will restrict to maximal building sets, that is, $\G = \L_{\Harr} \setminus \{V\}$. This slightly differs from \cite{Axelrod_Singer_1993} and \cite{lambrechts_formality_2012} (the braid case) where the authors consider the minimal building set instead. This does have the effect of changing the structure of the strata of $\SFM[\Br_n]$, which we will come back to later. \\

For now let us go back to the example of braid arrangements. For any $m \geq n$ we have a fibration $\pi^{m}_{n}:\Conf_m(\Cx) \twoheadrightarrow \Conf_n(\Cx)$ given by forgetting the last $m-n$ points. This fibration induces a semi-algebraic oriented bundle $\overline{\pi}^{m}_{n}: \SFM[\Br_m]\twoheadrightarrow \SFM[\Br_n]$ at the level of compactifications (see \cite[Section 8]{Hardt2008RealHT} for the definition of a semi-algebraic oriented bundle), which induces the pushforward map $(\overline{\pi}^m_n)_{\star}: \Omega_{min}^{\bullet}(\SFM[\Br_m], \R) \rightarrow \Omega_{PA}^{\bullet - 2(m-n)}(\SFM[\Br_n], \R)$ given by integrating along the fibers of $\opi^m_n.$ Here $\Omega^{\bullet}_{min}(\SFM[\Br_m], \R)$ denotes the cdga of minimal semi-algebraic forms of $\SFM[\Br_n]$ with real coefficients, which should be thought of as linear combinations $\sum f_0 \d f_1\wedge\cdots \d f_n$ with semi-algebraic $f_i$'s. On the other hand, $\Omega_{PA}^{\bullet}(\SFM[\Br_n], \R)$ denotes the cdga of piecewise algebraic forms, which should be thought of as push-forwards of minimal forms. In this article we will treat $\Omega^{\bullet}_{min}$ and $\Omega^{\bullet}_{PA}$ as black boxes, referring to \cite{Hardt2008RealHT} for more details. We are now in position to describe the qiso of cdgas $\I: \Gra(n) \xrightarrow{\sim} \Omega_{PA}^{\bullet}(\SFM[\Br_n], \R)$ alluded to at the end of Section \ref{secgra}. This qiso is defined by setting for all diagran $\Gamma$ with $n$ external vertices and $m-n$ internal ones: $$\I(\Gamma) = (\overline{\pi}^m_n)_{\star}\left(\bigwedge_{(i,j) \textrm{ edge of } \Gamma} \omega_{ij}\right),$$ where $\omega_{ij}$ denotes the form $\frac{\d(z_i - z_j)}{z_i - z_j}$ and the internal vertices of $\Gamma$ are labeled from $m-n+1$ to $m$ in any order. The main take away is that $\Gamma \in \Gra(n)$ models an explicit form initially coming from a larger braid arrangement $\Br_m$ and then partially integrated back down to $\Br_n$. Let us try to reformulate this situation in the language of hyperplane arrangements. \\

For any flat $F$ of a hyperplane arrangement $\Harr$, the subset of hyperplanes of $\Harr$ containing $F$ induces a hyperplane arrangement in the quotient $V/F$, denoted $\Harr_{F}$ and called the \textit{restriction} of $\Harr$ at $F$. The quotient map $V \twoheadrightarrow V/F$ induces a restriction map $\pi_F: \A_{\Harr}\rightarrow \A_{\Harr_F}.$ For instance if $\Harr$ is the braid arrangement $\Br_m$ and $F$ is the flat $\{z_1= \cdots = z_n\}$ for some $n\leq m$ then $\Harr_{F}$ can be identified with $\Br_n$ and under this identification the map $\A_{\Harr} \rightarrow \A_{\Harr_{F}}$ is the map forgetting the last $m-n$ points. In this particular case that map is a fibration but this is not true in general. Even for braid arrangements if we consider the flat $F =\{z_1=z_2\}\cap \{z_3=z_4\}$ of $\Br_4$ then the corresponding restriction map is not a fibration. A classical fact, which will be central for us in this article, is that detecting whether $\pi_F$ is a fibration can be done at the combinatorial level, that is, simply by analyzing the properties of $F$ as an element of the intersection lattice $\L_{\Harr}$. For this we need the notion of modularity. 

\begin{madef}[Modular element]
Let $F$ be some element in some intersection lattice $\L_{\Harr}$. We say that $F$ is \textit{modular} if for every $F'\in \L_{\Harr}$ we have the equality 
\begin{equation*}
    \codim(F) + \codim(F') = \codim (F\wedge F') + \codim (F\vee F'),
\end{equation*}
where $F\wedge F'$ and $F \vee F'$ denote respectively the supremum and infimum of $F$ and $F'$ in the lattice $\L_\Harr$.
\end{madef}
Beware that while $F\vee F'$ is simply the intersection of $F$ and $F'$ (which is again a flat of $\Harr$), the infimum $F\wedge F'$ is not in general $F+F'$ (which has a priori no reason to be a flat of $\Harr$) but only the biggest element of $\L_{\Harr}$ containing $F+F'$. The following result is due to Paris \cite{Paris_2000}, generalizing an earlier result of Terao \cite{Terao_1986} in the dimension $1$ case. See also Falk--Proudfoot \cite{Falk_2002}. 
\begin{theo}[\cite{Paris_2000}]
For any complex hyperplane arrangement $\Harr$ and any flat $F\in \L_{\Harr}$, the restriction map $\pi_F:\A_{\Harr}\rightarrow \A_{\Harr_F}$ is a fibration if and only if $F$ is a modular element of $\L_{\Harr}$. 
\end{theo}
For instance, going back to $F=\{z_1= z_2\}\cap \{z_3=z_4\}$ in $\L_{\Br_4}$, let us consider the flat $F'= \{z_1 = z_3\}\cap \{z_2= z_4\}$. We have $F\vee F' = \{0\}$ (in the essentialization) and we have $F\wedge F' = \Cx^4/\{z_1=z_2=z_3=z_4\}$ the full vector space which gives 
\begin{align*}
    \codim(F\vee F') + \codim(F\wedge F')  = 3 \neq 4 = \codim(F) + \codim(F'). 
\end{align*}
This shows that $F$ is not modular which explains why $\pi_F$ is not a fibration. If $F$ is a modular flat of a complexified hyperplane arrangement $\Harr$, the fibration $\pi_F$ induces a semi-algebraic oriented bundle $\overline{\pi}_F:\SFM[\Harr]\twoheadrightarrow \SFM[\Harr_{F}]$, with fiber of real dimension $2\dim(F)$. We shall call those semi-algebraic oriented bundles ``modular projections''. Assume that for each hyperplane $H$ in a complexified hyperplane arrangement $\Harr$ inside a vector space $V$, we have chosen a linear form $f_H:V\rightarrow \Cx$ such that $\ker f_H = H$. This gives explicit isomorphisms $\widetilde{f_H}:\S(V/H) \xrightarrow{\sim} \S^1 $ for all $H \in \Harr$. Besides, by definition of $\SFM[\Harr]$ for all linearly ordered multiset $J\subset \Harr$ we have a morphism $\theta_J: \SFM[\Harr]\rightarrow \prod_{H\in J}\S(V/H)$. If we let $\vol_J$ denote the standard volume form on $(\S^1)^J$, we can set $\omega_{\Harr, J} := \theta_J^\star(\prod_{H\in J}\widetilde{f_H})^*\vol_J$ which is a minimal semi-algebraic form of $\SFM[\Harr]$. Assume that we have fixed a complexified hyperplane arrangement $\Harr$ in a complex vector space $V$. Assume also that we are given a second complexified hyperplane arrangement $\bigarr$ in a vector space $V'$, together with a modular flat $F$ of $\bigarr$ and an isomorphism of hyperplane arrangements $\iota:\Harr\xrightarrow{\sim} \bigarr/F$ induced by a linear isomorphism $\iota: V \xrightarrow{\sim} V'/F$. For any linearly ordered multi-set $J\subset \bigarr$ we can consider the semi-algebraic form 
\begin{equation*}
    \omega_{\bigarr,\iota,F,J} \coloneqq \iota^{\star}(\overline{\pi}_F)_{\star}\omega_{\bigarr,J} \in \Omega_{PA}^{|J|-2\dim F}(\SFM[\Harr], \R).
\end{equation*}
The rest of this section is devoted to highlighting the main properties of the forms $\omega_{\bigarr,\iota, F, J}$, which will then be used as definition of our combinatorial model in the next section. All the proofs are direct adaptations of the proofs given in \cite[Section 9]{lambrechts_formality_2012}. We start with the following lemma.

\begin{lemma}\label{lemmarelform}
Let $\bigarr, \iota, F, J$ be as above.  If there exists a proper modular flat $F'$ of $\bigarr$ contained in $F$ and such that $J$ has strictly less than $2$ hyperplanes not containing $F'$, then $\omega_{\bigarr, \iota, F, J}$ is zero. 
\end{lemma}
\begin{proof}
To simplify notations let us assume that $\iota$ is the identity. One can easily check that the flat $F$ is also modular in $\bigarr_{F'}$, which gives a modular projection $\overline{\pi}_{F',F}: \SFM[\bigarr_{F'}] \twoheadrightarrow \SFM[\bigarr_{F}]$. 
If $J$ contains no hyperplanes not containing $F'$, then we have the following commutative diagram of semi-algebraic maps: 
\begin{equation*}
    \begin{tikzpicture}
        \node[] (A) at (0,0) {$\SFM[\bigarr]$};
        \node[] (B) at (4,0) {$\SFM[\bigarr_{F'}]$};
        \node[] (C) at (8,0) {$(\S^1)^I$};
        \node[] (D) at (6,2) {$\SFM[\bigarr_{F}]$};
        \draw[->] (A) -- (B) node[midway,below, scale = 0.8] {$\overline{\pi}_{F'}$}; 
        \draw[->] (B) -- (C) node[midway,below, scale = 0.8] {$\prod_{H \in I}\widetilde{f}_{H}$}; 
        \draw[->] (A) -- (D) node[midway,above, scale = 0.8] {$\overline{\pi}_F$}; 
        \draw[->] (B) -- (D) node[pos=0.4,right= 3pt, scale = 0.8] {$\overline{\pi}_{F',F}$}; 
    \end{tikzpicture}
\end{equation*}
This means that we have 
\begin{equation*}
    \omega_{\bigarr, \Id, F,J} = (\overline{\pi}_F)_{\star}(\overline{\pi}_{F'})^{\star}(\prod_{H \in J}\widetilde{f}_{H})^{\star}(\vol_J).
\end{equation*}
For all $x\in \SFM[\bigarr_{F}]$ we have $\dim \opi_{F',F}^{-1}(x) < \dim \opi_F^{-1}(x)$, and so by \cite[Proposition 8.14]{Hardt2008RealHT} the morphism $(\overline{\pi}_F)_{\star}(\overline{\pi}_{F'})^{\star}$ is zero on minimal semi-algebraic forms. 

If $J$ contains a unique hyperplane $H$ not containing $F'$, we have the following commutative diagram of semi-algebraic maps: 
\begin{equation*}
    \begin{tikzpicture}
        \node[] (A) at (0,0) {$\SFM[\bigarr]$};
        \node[] (B) at (4,0) {$\S^1\times \SFM[\bigarr_{F'}]$};
        \node[] (C) at (8,0) {$(\S^1)^I$};
        \node[] (D) at (6,2) {$\SFM[\bigarr_{F}]$};
        \draw[->] (A) -- (B) node[midway,below, scale = 0.8] {$\widetilde{f}_H\times \overline{\pi}_{F'}$}; 
        \draw[->] (B) -- (C) node[midway,below, scale = 0.8] {$\Id\times \prod_{H' \in J\setminus H}\widetilde{f}_{H'}$}; 
        \draw[->] (A) -- (D) node[midway,above, scale = 0.8] {$\overline{\pi}_F$}; 
        \draw[->] (B) -- (D) node[pos=0.4,right, scale = 0.8] {$\overline{\pi}_{F',F}\circ \pr_2$}; 
    \end{tikzpicture}
\end{equation*}
This means that we have 
\begin{equation*}
    \omega_{\bigarr, \Id, F,J} = (\overline{\pi}_F)_{\star}(\widetilde{f}_H\times \overline{\pi}_{F'})^{\star}(\Id\times \prod_{H' \in J\setminus H}\widetilde{f}_{H'})^{\star}(\vol_J).
\end{equation*}
For all $x\in \SFM[\bigarr_{F}]$ we have $\dim (\opi_{F',F}\circ \pr_2)^{-1}(x) < \dim \opi_F^{-1}(x)$, and so by \cite[Proposition 8.14]{Hardt2008RealHT} the morphism $(\overline{\pi}_F)_{\star}(\widetilde{f}_H\times \overline{\pi}_{F'})^{\star}$ is zero on minimal semi-algebraic forms.
\end{proof}
This lemma motivates getting rid of diagrams with isolated vertices and univalent internal vertices in $\Gra(n)$. Killing diagrams with bivalent internal vertices is motivated by Lemma \cite[Lemma 9.11]{lambrechts_formality_2012}. We were not able to prove a generalization of that statement in the setting of hyperplane arrangements so we leave it as a conjecture. 
\begin{conj}
Let $\bigarr, \iota, F, J$ be as above.  If there exists a proper modular flat $F'$ of $\bigarr$ contained in $F$ and such that $J$ contains exactly $2$ hyperplanes not containing $F'$, then $\omega_{\bigarr, \iota, F, J}$ is zero. 
\end{conj}
We also have the following lemma. 
\begin{lemma}\label{lemmagradzero}
Let $\bigarr, \iota, F, J$ be as above. If there exists a modular flat $F'$ of $\bigarr$ such that $F' \wedge F = V$, $F'\vee F \neq \{\zero\}$ and such that every hyperplane of $J$ contains $F'$, then $\omega_{\bigarr, \iota, F,J}$ is zero. 
\end{lemma}
\begin{proof}
In that case we have the following commutative diagram of semi-algebraic maps: 
\begin{equation*}
    \begin{tikzpicture}
        \node[] (A) at (0,0) {$\SFM[\bigarr]$};
        \node[] (B) at (4,0) {$\SFM[\bigarr_{F}]\times \SFM[\bigarr_{F'}]$};
        \node[] (C) at (8,0) {$(\S^1)^J$};
        \node[] (D) at (6,2) {$\SFM[\bigarr_{F}]$};
        \draw[->] (A) -- (B) node[midway,below, scale = 0.8] {$\overline{\pi}_F\times \overline{\pi}_{F'}$}; 
        \draw[->] (B) -- (C) node[midway,below, scale = 0.8] {$(\prod_{H \in J}\widetilde{f}_{H})\circ \pr_2$}; 
        \draw[->] (A) -- (D) node[midway,above, scale = 0.8] {$\overline{\pi}_F$}; 
        \draw[->] (B) -- (D) node[pos=0.4,right, scale = 0.8] {$\pr_1$}; 
    \end{tikzpicture}
\end{equation*}
This means that we have 
\begin{equation*}
    \omega_{\bigarr, \Id, F, J} = (\overline{\pi}_F)_{\star}(\overline{\pi}_{F}\times \overline{\pi}_{F'})^{\star}((\prod_{H \in J}\widetilde{f}_{H})\circ \pr_2)^{\star}(\vol_J).
\end{equation*}
For all $x\in \SFM[\bigarr_{F}]$ we have 
$$\dim \pr_1^{-1}(x) = \dim \SFM[\bigarr_{F'}] = \codim F' < \dim F = \dim \opi_F^{-1}(x),$$ 
and so by \cite[Proposition 8.14]{Hardt2008RealHT} the morphism $(\overline{\pi}_F)_{\star}(\overline{\pi}_{F}\times \overline{\pi}_{F'})^{\star}$ is zero on minimal semi-algebraic forms. 
\end{proof}
We are now interested in describing how the forms $\omega_{\bigarr, \iota, F,J}$ multiply in $\Omega^{\bullet}_{PA}(\SFM[\Harr], \R).$ Suppose that we have two complexified hyperplane arrangements $\bigarr_1 = \{H^1_1, \ldots, H^1_{k_1}\}$ and $\bigarr_2= \{H^2_1, \ldots ,H^2_{k_2}\}$ in $V_1$ and $V_2$ respectively, together with modular flats $F_1$ and $F_2$ of $\bigarr_1$ and $\bigarr_2$ respectively, such that we have two isomorphisms of hyperplane arrangements $\iota_1:\Harr \xrightarrow{\sim} \bigarr_{F_1}$, $\iota_2:\Harr \xrightarrow{\sim} \bigarr_{F_2}$. Inside the pull-back vector space $V_3 = V_1\times_{\iota_1,\iota_2} V_2$ we have hyperplanes of the form $\{(v_1,v_2) \in V_3 \, | \, v_1 \in H \}$ where $H$ is some hyperplane of $\bigarr_1$ not containing $F_1$. Similarly we also have hyperplanes of the form $\{(v_1,v_2) \in V_3 \, | \, v_2 \in H\}$ for some hyperplane $H$ of $\bigarr_2$ not containing $F_2$ and finally we also have hyperplanes of the form $\{(v_1,v_2) \in V_3\, | \, v_1\in H, v_2 \in \overline{\imath_2\circ \imath_1^{-1}}H\}$ where $H$ is some hyperplane of $\bigarr_1$ containing $F$ and $\overline{\imath_2\circ \imath_1^{-1}}H$ is the unique hyperplane of $\bigarr_2$ such that we have $\imath_2^{-2}(\overline{\imath_2\circ \imath_1^{-1}}H/F_2) = \imath_1^{-1}(H/F_1)$. Let us denote by $\bigarr_1\cup_{\iota_1,\iota_2}\bigarr_2$ the arrangement given by those three types of hyperplanes in the vector space $V_3$. One can check that the flat $F_{12} = \{(v_1,v_2) \in V_3 \, | \, v_1 \in F_1, v_2 \in F_2\}$ is modular in $\bigarr_1\cup_{\imath_1\imath_2} \bigarr_2$ and we have an isomorphism of hyperplane arrangements $\imath_{12}: \Harr \simeq (\bigarr_1\cup_{\imath_1\imath_2} \bigarr_2)_{F_{12}}$. Note that as a set $\bigarr_1\cup_{\imath_1\imath_2}\bigarr_2$ is simply the pushout of $\bigarr_1$ and $\bigarr_2$ along $\Harr$. Let us denote $F'_1 = \{(v_1, v_2) \in V_3\, | \, v_1 = 0\},$ and $F'_2 = \{(v_1, v_2) \in V_3\, | \, v_2 = 0\}$. For $i=1,2$ the subspace $F'_i$ is a modular flat of $\bigarr_1\cup_{\iota_1,\iota_2}\bigarr_2$ and we have an isomorphism of hyperplane arrangements $\imath'_i:\bigarr_i \simeq (\bigarr_1\cup_{\iota_1,\iota_2}\bigarr_2)_{F'_i}$. For any linearly ordered multisets $J, J'$ of hyperplanes of $\bigarr_1$ and $\bigarr_2$ respectively, let us denote by $J_1\bullet J_2$ the multiset union $J_1\cup J_2 \subset \bigarr_1\cup_{\imath_1\imath_2}\bigarr_2$ linearly ordered by concatenating the orders on $J_1$ and $J_2$, putting the elements of $J_1$ before the elements of $J_2$. We have the following lemma. 
\begin{lemma}
Let $(\bigarr_1, \imath_1, F_1)$ and $(\bigarr_2, \imath_2,F_2)$ be as in the preceding paragraph. For any linearly ordered multisets $J_1, J_2$ of hyperplanes of $\Harr_1$ and $\Harr_2$ respectively, we have the equality of forms 
\begin{equation*}
    \omega_{\bigarr_1,\imath_1,F_1, J_1}\wedge \omega_{\bigarr_2, \imath_2, F_2, J_2} = \omega_{\bigarr_1\cup_{\imath_1\imath_2}\bigarr_2, \imath_{12}, F_{12}, J_1\bullet J_2} 
\end{equation*}
in $\Omega_{PA}^{\bullet}(\SFM[\Harr], \R)$. 
\end{lemma}
\begin{proof}
Notice that we have a pullback diagram 
\begin{equation*}
    \begin{tikzpicture}
        \node[] (A) at (0,0) {$\SFM[\bigarr_1\cup_{\imath_{12}}\bigarr_2]$};
        \node[] (B) at (5,0) {$\SFM[\bigarr_1]$};
        \node[] (C) at (0,-3) {$\SFM[\bigarr_2]$};
        \node[] (D) at (5,-3) {$\SFM[\Harr]$};
        \draw[->] (A) -- (B) node[midway, above] {$\imath_1^{,-1}\overline{\pi}_{F'_1}$};
        \draw[->] (A) -- (C) node[midway, left] {$\imath^{,-1}_2\overline{\pi}_{F'_2}$};
        \draw[->] (B) -- (D) node[midway, right] {$\imath_1^{-1}\overline{\pi}_{F_1}$};
        \draw[->] (C) -- (D) node[midway, below] {$\imath_2^{-1}\overline{\pi}_{F_2}$};
        
    \end{tikzpicture}
\end{equation*}
This implies that we have the equalities 
\begin{align*}
    \omega_{\bigarr_1\cup_{\imath_1\imath_2}\bigarr_2, \imath_{12}, F_{12}, J_1\bullet J_2} &= \iota_{12}^{\star}(\overline{\pi}_{F_{12}})_{\star}(\theta_{J_1\bullet J_2})^{\star}(\prod_{H \in J_1\bullet J_2}\widetilde{f}_H)^{\star}(\vol_{J_1\bullet J_2}) \\
    &= \iota_{12}^{\star}(\overline{\pi}_{F_{12}})_{\star}\left(\biggl(\overline{\pi}_{F'_1}^{\star}\theta_{J_1}^{\star}(\prod_{H \in J_1} \widetilde{f}_H)^{\star}\vol_{J_1}\biggr) \wedge \biggl( \overline{\pi}_{F'_2}^{\star}\theta_{J_2}^{\star}(\prod_{H \in J_2} \widetilde{f}_H)^{\star}\vol_{J_2}\biggr)\right) \\ 
    &= \biggl(\imath_{12}^{\star}(\imath_1^{-1}\overline{\pi}_{F_1})_{\star}\theta_{J_1}^{\star}(\prod_{H \in J_1} \widetilde{f}_H)^{\star}\vol_{J_1} \biggr)\wedge \biggl(\imath_{12}^{\star}(\imath_2^{-1}\overline{\pi}_{F_2})_{\star}\theta_{J_2}^{\star}(\prod_{H \in J_2} \widetilde{f}_H)^{\star}\vol_{J_2}\biggr) \\
    &= \omega_{\bigarr_1,\imath_1,F_1, J_1}\wedge \omega_{\bigarr_2, \imath_2, F_2, J_2},
\end{align*}
where the third equality comes from \cite[Proposition 8.15]{Hardt2008RealHT}. 
\end{proof}
Notice that when we restrict to braid arrangements the proof of the analogous statement \cite[Proposition 9.2]{lambrechts_formality_2012} is slightly more complicated because the pushout arrangement $\Br_p\cup_{\Br_{n}}\Br_q$ is not a braid arrangement, it is the graphical arrangement associated to the graph given by gluing the complete graphs $K_p, K_q$ along a subcomplete graph $K_n$. The spherical compactification $\SFM[\Br_p\cup_{\Br_{n}}\Br_q]$ is the singular configuration space described in \cite[Section 5.5]{lambrechts_formality_2012}. In that setting one has to relate the singular configuration space to a bigger braid arrangement, which is not necessary for us because we are considering all arrangements at once.  \\

Our final task is to compute the differential of the forms $\omega_{\bigarr, \imath, F, J}$ in the differential complex $\Omega^{\bullet}_{PA}(\SFM[\Harr], \R).$ If $F$ is a flat of a hyperplane arrangement $\bigarr$ in some vector space $V$, then one can consider the set of subspaces $\bigarr^F \coloneqq \{F\cap H \, | \,  H \in \bigarr, H\nsupseteq F \}$, which is a hyperplane arrangement in $F$ called the \textit{contraction of $\bigarr$ at $F$}. If in addition we have a modular flat $F'$ of $\bigarr$ such that $F +  F' = V$ then we have an isomorphism of hyperplane arrangements $\iota_{F,F'}: \bigarr_{F'} \simeq \bigarr_{F\vee F'}^{F}$ induced by the linear isomorphism $\frac{V}{F'} \simeq \frac{F}{F\cap F'}.$ We have the following lemma, which generalizes \cite[Proposition 9.12]{lambrechts_formality_2012}. 
\begin{lemma}
Let $\bigarr$ be a complexified hyperplane arrangement, $F$ a modular flat of $\bigarr$ and $\iota: \Harr \xrightarrow{\sim} \bigarr_{F}$ an isomorphism of hyperplane arrangements. For all linearly ordered multiset $J = \{H_1 \vartriangleleft \cdots \vartriangleleft H_k\} \subset \bigarr$ we have the equality
\begin{equation*}
    \d(\omega_{\bigarr, \iota, F, J}) = \sum_{\substack{1 \leq i \leq k \\ H_i \nsupseteq F }}(-1)^{i+1}\omega_{\bigarr^{H_i}, \iota_{H_i,F}\circ\iota, F\vee H_i, H_i\vee (J\setminus H_i)},
\end{equation*}
where $H_i \vee (H\setminus H_i)$ denotes the linearly ordered multiset $\{H_i\vee H_1 \vartriangleleft \cdots \vartriangleleft \widehat{H_i\vee H_i} \vartriangleleft \cdots \vartriangleleft H_i\vee H_k\}.$
\end{lemma}
The proof is slightly technical and is a direct adaptation of that given in \cite{lambrechts_formality_2012} so we streamline it.  
\begin{proof}[Sketch of proof]
To simplify notations we assume $\iota = \Id$. Denote by $\overline{\pi}_F^{\partial}:\SFM^{\partial,F}[\bigarr] \twoheadrightarrow \SFM[\bigarr_F]$ the fiberwise boundary of the fibration $\overline{\pi}_F: \SFM[\bigarr] \twoheadrightarrow \SFM[\bigarr_F]$. In short $\SFM^{\partial,F}[\bigarr]$ is the union of the boundaries of the fibers of $\overline{\pi}_F$. We refer to \cite[Section 8]{Hardt2008RealHT} for more details on this construction. By the fiberwise Stokes formula \cite[Proposition 8.12]{Hardt2008RealHT} we have 
\begin{align*}
    \d(\omega_{\bigarr, \Id, F, J}) &= \d\biggl((\overline{\pi}_F)_{\star}\theta_J^{\star}(\prod_{H\in J}\widetilde{f}_H)^{\star}(\vol_J)\biggr) \\ 
    &= (-1)^{\deg(\omega_{\bigarr, \Id, F, J})} (\overline{\pi}_F^{\partial})_{\star}\Bigl(\theta_J^{\star}(\prod_{H\in J}\widetilde{f}_H)^{\star}(\vol_J)\, \Big|  \, \SFM^{\partial,F}[\bigarr]\Bigr)
\end{align*}
where $\bullet \, \big|\, \SFM^{\partial,F}[\bigarr]$ denotes the restriction to $\SFM^{\partial,F}[\bigarr]$. Gaiffi \cite{Gaiffi_2003} proved that the boundary $\partial \SFM[\Harr]$ is a union of codim $2$ submanifolds $\partial_G\SFM[\bigarr]$ indexed by proper flats of $\bigarr$, and those submanifolds have transverse intersection. Moreover, we have isomorphisms $\Phi_{G}: \SFM[\bigarr_{G}]\times \SFM[\bigarr^G]\xrightarrow{\sim} \partial_G\SFM[\bigarr] $. More generally if we consider the spherical compactification associated to a building set $\G$ the codim $2$ stratas are indexed by proper flats of $\G.$ The fiberwise boundary $\SFM^{\partial, F}$ is the union of those codim $2$ submanifolds with $G$ running over the subset $\mathcal{BF}(F)$ of flats $G$ such that we have either $G \subset F$ or $F+G = V.$ For such a $G$ denote $$\omega^{\partial}_G \coloneqq \theta_J^{\star}(\prod_{H\in J}\widetilde{f}_H)^{\star}(\vol_J)\, \Big|  \, \partial_G\SFM[\bigarr].$$ By additivity of integration along the fiber \cite[Proposition 8.11]{Hardt2008RealHT} we have the equality 
\begin{equation*}
    (\overline{\pi}_F^{\partial})_{\star}\Bigl(\theta_J^{\star}(\prod_{H\in J}\widetilde{f}_H)^{\star}(\vol_J)\, \Big|  \, \SFM^{\partial,F}[\bigarr]\Bigr) = \sum_{G\in \mathcal{BF}(F)}(\overline{\pi}^{\partial}| \partial_G\SFM[\bigarr]) _{\star}\omega^{\partial}_G.
\end{equation*}
Up to some sign, for all $G \in \mathcal{BF}(F)$ we have the equality 
\begin{equation}\label{eqomegapartial}
     (\overline{\pi}^{\partial}| \partial_G\SFM[\bigarr]) _{\star}\omega^{\partial}_G = (\overline{\pi}^{\partial}\circ \Phi_G)_{\star}(\omega'_G\times \overline{\omega}_G)
\end{equation} 
for some forms $\omega'_G \in \Omega_{PA}(\SFM[\bigarr_G], \R)$ and $\overline{\omega}'_G \in \Omega_{PA}(\SFM[\bigarr^G], \R).$ For $G\in \mathcal{BF}(F)$ such that $G \subset F$ we have the equality 

\begin{equation}\label{eqpipartial}
    \overline{\pi}^{\partial}| \partial_G\SFM[\bigarr] = \overline{\pi}_{F}^{G}  \circ\pr_1\circ\Phi_{G}^{-1},
\end{equation} 
where $\pr_1: \SFM[\bigarr_G]\times \SFM[\bigarr^G]$ is the projection on the first factor. Equations \eqref{eqomegapartial} and \eqref{eqpipartial} give
\begin{equation*}
     (\overline{\pi}^{\partial}| \partial_G\SFM[\bigarr]) _{\star}\omega^{\partial}_G = (\overline{\pi}^G_F)_{\star}\omega'_G\cdot\langle \SFM[\bigarr^G],\overline{\omega}_G\rangle.
\end{equation*}
One can then prove that we have in fact $\langle \SFM[\bigarr^G],\overline{\omega}_G\rangle = 0.$ On the other hand if $G\in \mathcal{BF}(F)$ is such that we have $G\wedge F = V$, using similar arguments one can prove that the term $(\overline{\pi}^{\partial}| \partial_G\SFM[\bigarr]) _{\star}\omega^{\partial}_G$ is zero unless $G$ is some hyperplane $H$ in $J$ not containing $F$, in which case we get $$(\overline{\pi}^{\partial}| \partial_G\SFM[\bigarr]) _{\star}\omega^{\partial}_G = \pm \omega_{\bigarr^{H}, \iota_{H,F}\circ\iota, F\vee H, H\vee (J\setminus H)},$$
which concludes the proof.
\end{proof}
\section{The combinatorics}\label{secMD}
In this section we construct a quasi-free cdga analogous to Kontsevich's cdga of admissible graphs and inspired by the results of the previous section, in a purely combinatorial setting. This means that we will not be dealing directly with hyperplane arrangements but with a certain class of posets modeling hyperplane arrangements. We start with some combinatorial preliminaries. 

\begin{madef}[Lattice]
    A poset $\L$ is called a \textit{lattice} if every pair of elements $F_1, F_2$ admits a supremum, denoted $F_1\vee F_2$, and an infimum, denoted $F_1 \wedge F_2$. 
\end{madef}
This implies more generally that any finite subset $S \subset \L$ admits a supremum and an infimum, which will be denoted by $\medvee S$ and $\bigwedge S$ respectively. If $\L$ is finite this implies that $\L$ is bounded, with upper bound $\bigvee \L $ denoted $\un$ and lower bound $\bigwedge \L$ denoted $\zero$. We shall be interested in the following class of lattices, which mimic intersection lattices of hyperplane arrangements. 
\begin{madef}[Geometric lattice]
    A finite lattice $\L$ is called \textit{geometric} if it satisfies the following properties: 
    \begin{enumerate}[label={\roman*)}]
        \item $\L$ is well-ranked, meaning that for every $F$ in $\L$, every maximal chain from $\zero$ to $F$ has the same cardinality, denoted $\rk F$.
        \item $\L$ is atomic, meaning that every element $F\in \L$ can be obtained as a join of rank $1$ elements (also called atoms).
        \item $\L$ is semi-modular, meaning that for all $F_1,F_2\in \L$ we have the inequality 
        \begin{equation*}
            \rk(F_1\wedge F_2) + \rk(F_1\vee F_2) \leq \rk(F_1) + \rk(F_2).
        \end{equation*}
    \end{enumerate}
\end{madef}
A geometric lattice is called \textit{non-trivial} if it is neither empty nor a singleton. The set of atoms of a geometric lattice $\L$ will be denoted by $\At(\L)$. We extend the rank function to sets of atoms by setting $\rk  \coloneqq \rk \bigvee J$ for all $J \subset \At(\L).$ The rank of $\L$ is simply $\rk \L \coloneqq \rk \un$. An element of $\L$ will be called a \textit{flat}. A flat will be called \textit{proper} if it is neither $\zero$ nor $\un$. A \textit{coatom} of $\L$ is a flat of rank $\rk\L - 1$. For any $F$ in $\L$ we will also use $F$ to denote the set of atoms $\{H \in \At(\L) \, | \, H \leq F \}$ (leaving to  the context to clear up any possible confusion). 
\begin{ex}\label{exgeolatt}
Consider the following set of points and lines. 
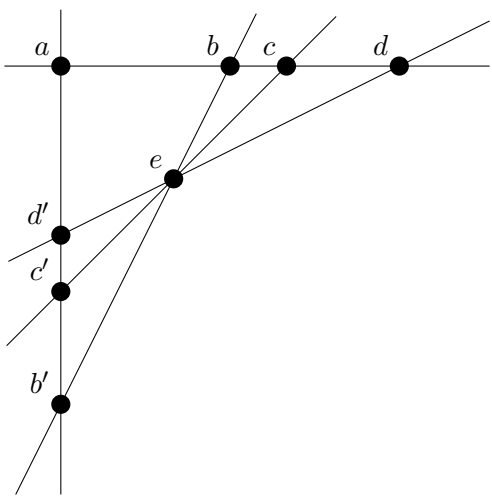
\begin{equation*}
\begin{tikzpicture}
    \node[circle, fill, black, inner xsep=0.9mm,inner ysep = 0.9mm]  (A) at (0,0) {} node[above left =-3pt of A, outer sep=0pt] {$a$};
    \node[circle, fill, black, inner xsep=0.9mm,inner ysep = 0.9mm]  (B) at (3,0) {};
    \node[circle, fill, black, inner xsep=0.9mm,inner ysep = 0.9mm]  (C) at (0,-3) {};
    \node[circle, fill, black, inner xsep=0.9mm,inner ysep = 0.9mm]  (D) at (1.5,-1.5) {};\node[circle, fill, black, inner xsep=0.9mm,inner ysep = 0.9mm]  (E) at (2.25,0) {};
    \node[circle, fill, black, inner xsep=0.9mm,inner ysep = 0.9mm]  (F) at (0,-2.25) {};
    \node[circle, fill, black, inner xsep=0.9mm,inner ysep = 0.9mm]  (I) at (4.5,0) {};
    \node[circle, fill, black, inner xsep=0.9mm,inner ysep = 0.9mm]  (J) at (0,-4.5) {};
    \node[above left =-3pt of E, outer sep=0pt] {$b$};
    \node[above left =-3pt of B, outer sep=0pt] {$c$};
    \node[above left =-3pt of I, outer sep=0pt] {$d$};
    \node[above left =-3pt of F, outer sep=0pt] {$d'$};
    \node[above left =-3pt of C, outer sep=0pt] {$c'$};
    \node[above left =-3pt of J, outer sep=0pt] {$b'$};
    \node[above left =-3pt of D, outer sep=0pt] {$e$};
    \draw[black] ($ (A)!-0.25!(B) $) -- ($ (A)!1.9!(B) $); 
    \draw[black] ($ (A)!-0.25!(C) $) -- ($ (A)!1.9!(C) $);
    \draw[black] ($ (B)!-0.22!(C) $) -- ($ (B)!1.24!(C) $);
    \draw[black] ($ (E)!-0.465!(D) $) -- ($ (E)!3.8!(D) $);
    \draw[black] ($ (F)!-0.465!(D) $) -- ($ (F)!3.8!(D) $);
\end{tikzpicture}
\end{equation*}
This picture defines a geometric lattice of rank $3$, with atoms given by points and coatoms given by the drawn lines and all two point lines. The order is given by the incidence relations. 
\end{ex}
\begin{ex}\label{exgeolattgraph}
    If $G$ is a simple graph with set of edges $E$ and set of vertices $V$, we can associate to $G$ a geometric lattice $\L_G$ defined as the subposet of the power set of $E$ (ordered by inclusion) consisting of subsets $F\subset E$ satisfying the property that for all edges $(i,j) \in E$, if $F$ contains a set of edges forming a path from $i$ to $j$ then $(i,j)$ belongs to $F$. This poset is also the intersection lattice of the graphical hyperplane arrangement associated to $G$.
\end{ex}
It is easy to check that for any hyperplane arrangement $\Harr$ (not necessarily complex), the intersection lattice $\L_{\Harr}$ is a geometric lattice. The atoms of $\L_{\Harr}$ are the hyperplanes of the arrangement. The lattice $\L_{\Harr}$ is well-ranked by the existence of dimension in linear algebra, it is atomic by definition and it is semi-modular by the classical formula 
$$ \dim (F_1\cap F_2) + \dim(F_1 + F_2)  = \dim(F_1) + \dim(F_2).$$
However, not all geometric lattices are intersection lattices, which means that we won't be able to use any result from the previous section. However, those results will still be our main source of inspiration. Geometric lattices are equivalent to simple matroids (see \cite{oxley} for a general reference), which have many other axiomatizations. In this article we chose the axiomatization by geometric lattices because it is the closest to hyperplane arrangements. Since we do not have any arrangement anymore we cannot consider the arrangement complement but we can still associate to a geometric lattice $\L$ a graded commutative algebra mimicking the cohomology of the arrangement complement, called the Orlik--Solomon algebra of $\L$. For a geometric lattice $\L$, a subset $C\subset \At(\L)$ is called a \textit{circuit} if we have $\rk (C) = \#C -1 $ and $\rk(C\setminus c) = \#C -1$ for all $c \in C$. 
\begin{madef}[Orlik--Solomon algebra]
Let $\L$ be a geometric lattice. The \textit{Orlik--Solomon algebra} of $\L$, denoted $\OS(\L)$, is the graded commutative algebra over $\Q$ defined by the presentation 
\begin{equation*}
    \OS(\L) \coloneqq \frac{\bigwedge^{\bullet}[ e_H \, \, | \, \,  H \in \At(\L), \deg(e_H) = 1]}{(\delta e_C, C \textrm{ circuit of } \L) },
\end{equation*}
where $\delta$ is the only derivation of $\bigwedge^{\bullet}[e_H \, \, | \, \,  H \in \At(\L), \deg(e_H) = 1]$ sending $e_H$ to $1$ for all $H$, and $e_C$ means $\pm \prod_{H\in C} e_H$. 
\end{madef}
For instance for $\L = \L_{\Br_n}$, if for all $1\leq i<j\leq n $ we denote $H_{ij} = \{z_i = z_j\}$ viewed as an atom of $\L_{\Br_n}$, then for all $1 \leq i<j<k \leq n$ the set $\{H_{ij}, H_{ik}, H_{jk}\}$ is a circuit of $\L_{\Br_n}$, which gives the relation $e_{H_{ij}}e_{H_{ik}} - e_{H_{ij}}e_{H_{jk}} + e_{H_{ik}}e_{H_{jk}} = 0$ in $\OS(\L_{\Br_n})$. One can check that those relations generate all the others in $\OS(\L_{\Br_n})$ and so by equation \eqref{presArnold} we get $\OS(\L_{\Br_n}) \simeq H^{\bullet}(\Conf_n(\Cx), \Q).$ More generally we have the following seminal result due to Orlik and Solomon.
\begin{theo}[\cite{OS_1980}]
For all complex hyperplane arrangement $\Harr$, we have an isomomorphism of graded commutative algebras $\OS(\L_{\Harr}) \simeq \H^\bullet (\A_{\Harr}, \Q).$
\end{theo}
The element $e_H$ models the $1$-form $\d f_H/f_H$ for some chosen annihilator $f_H$ of $H$, which corresponds to ``turning around the hyperplane $H$''. The main goal of this section is to define a quasi-free dg model of $\OS(\L)$ for a subclass of geometric lattices called supersolvable, which we will introduce later. For now let us go back to recollections on geometric lattices. We have the following classical facts. 
\begin{prop}
For any geometric lattices $\L_1$ and $\L_2$ the product poset $\L_1\times \L_2$ is a geometric lattice. 
\end{prop}
A geometric lattice is called \textit{irreducible} if it is not a product of non-trivial geometric lattices.  
\begin{prop}\label{propdecir}
A geometric lattice decomposes uniquely as a product of irreducible geometric lattices. 
\end{prop}
\begin{prop}
For all $F_1, F_2$ in a geometric lattice $\L$, the subposet $[F_1, F_2] \coloneqq \{F \in \L \, | \,  F_1\leq F \leq F_2\} \subset \L$ is a geometric lattice. 
\end{prop}
When $F_1 = \zero$ we have a finer notion of inclusion given by the following definition. 
\begin{madeflemma}[Embedding]
Let $\varphi: \L_1 \rightarrow \L_2$ be an increasing map between two geometric lattices. We say that $\varphi$ is an \textit{embedding} if $\varphi$ is injective, compatible with the join on both side and sends atoms of $\L_1$ to atoms of $\L_2$. For any geometric lattice $\L$ and any subset $J\subset \At(\L)$ the poset $\L_{|J} \coloneqq \{ F \in \L \, | \, F \textrm{ can be obtained as a join of elements of } J \} \subset \L$ is a geometric lattice and the inclusion $\L_{|J}\subset \L$ is an embedding. That poset will be called the restriction of $\L$ at $J$. Any embedding factorizes as a composition of an isomorphism and a restriction.
\end{madeflemma}
\begin{ex}\label{exgeolattres}
The restriction of the geometric lattice $\L$ described in Example \ref{exgeolatt} to the subset of atoms $\{a,b,c,d,b',c',d'\}$ is the geometric lattice associated to the point-line configuration 
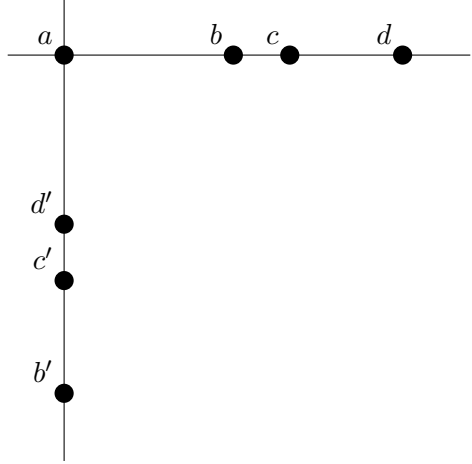
\begin{equation*}
    \begin{tikzpicture}
    \node[circle, fill, black, inner xsep=0.9mm,inner ysep = 0.9mm]  (A) at (0,0) {} node[above left =-3pt of A, outer sep=0pt] {$a$};
    \node[circle, fill, black, inner xsep=0.9mm,inner ysep = 0.9mm]  (B) at (3,0) {};
    \node[circle, fill, black, inner xsep=0.9mm,inner ysep = 0.9mm]  (C) at (0,-3) {};
    \node[circle, fill, black, inner xsep=0.9mm,inner ysep = 0.9mm]  (E) at (2.25,0) {};
    \node[circle, fill, black, inner xsep=0.9mm,inner ysep = 0.9mm]  (F) at (0,-2.25) {};
    \node[circle, fill, black, inner xsep=0.9mm,inner ysep = 0.9mm]  (I) at (4.5,0) {};
    \node[circle, fill, black, inner xsep=0.9mm,inner ysep = 0.9mm]  (J) at (0,-4.5) {};
    \node[above left =-3pt of E, outer sep=0pt] {$b$};
    \node[above left =-3pt of B, outer sep=0pt] {$c$};
    \node[above left =-3pt of I, outer sep=0pt] {$d$};
    \node[above left =-3pt of F, outer sep=0pt] {$d'$};
    \node[above left =-3pt of C, outer sep=0pt] {$c'$};
    \node[above left =-3pt of J, outer sep=0pt] {$b'$};
    \draw[black] ($ (A)!-0.25!(B) $) -- ($ (A)!1.8!(B) $); 
    \draw[black] ($ (A)!-0.25!(C) $) -- ($ (A)!1.8!(C) $);
\end{tikzpicture}
\end{equation*}
\end{ex}
\begin{ex}\label{exresgraph}
If $\L_G$ is the geometric lattice associated to some graph $G = (V,E)$ (see Example  \ref{exgeolattgraph}) and $J$ is some subset of $E$, the restriction $\L_{|J}$ is the geometric lattice associated to the restriction of the graph $G$ to $J.$
\end{ex}
\begin{madeflemma}[Modular element \cite{Stanley1971}]
    Let $F$ be an element in a geometric lattice $\L$. The following properties are equivalent: 
    \begin{enumerate}
        \item For all $F'\in \L$ we have the equality 
        \begin{equation*}
            \rk(F\wedge F') + \rk(F\vee F') = \rk(F) + \rk(F').
        \end{equation*}
        \item For all $A,B\in \L$ such that $A\leq F$, we have the equality $$F\wedge(A\vee B) = A \vee (F\wedge B).$$ 
        \item For all $A,B \in \L$ such that $A\leq B$, we have the equality $$ B\wedge(A\vee F) = A \vee (B\wedge F).$$
    \end{enumerate}
    If those properties are true we say that $F$ is \textit{modular}.
\end{madeflemma}
\begin{ex}\label{exmodular}
Note that an atom in a geometric lattice is always modular. At the other extreme, a coatom $F$ is modular if and only if taking the meet with $F$ drops the rank of a flat by at most one. For example for the geometric lattice $\L$ described in Example \ref{exgeolatt}, the line $\{a,b,c,d\}$ is modular in $\L$ because that line intersects every other line at a point of $\L$. On the other hand the line $\{a,e\}$ is not modular because that line has an empty intersection with the line $\{d',b\}$ for instance. 
\end{ex}
\begin{ex}\label{exmodulargraph}
Let $F$ be a flat of some geometric lattice $\L_G$ associated to graph $G= (V,E)$ (see Example \ref{exgeolattgraph}). Assume furthermore that $F$ is connected as a set of edges. By \cite[Proposition 4.3.2]{Brylawski_1975}, the flat $F$ is modular if and only if for all pairs of vertices $(u,v)$ in $V$, if $F$ contains a set of edges $J$ forming a path from $u$ to $v$ and there exists another path from $u$ to $v$ in $G$ which is vertex disjoint from $J$, then there exists an edge between $u$ and $v$ in $E$ (and so in $F$ as well). 
\end{ex}
We will need the following results about modularity.
\begin{lemma}[{\cite[Proposition 3.5]{Brylawski_1975}}]\label{lemmatransmod}
    Let $F_1 \leq F_2$ be two flats in a geometric lattice $\L$. If $F_1$ is modular in $[\zero,F_2]$ and $F_2$ is modular in $\L,$ then $F_1$ is modular in $\L$.
\end{lemma}
\begin{lemma}[{\cite[Proposition 3.8]{Brylawski_1975}}]\label{lemmaresmod}
Let $F$ be a modular flat of a geometric lattice $\L$. For any subset $J \subset \At(\L)$ containing $F$, $F$ is a modular flat of $\L_{|J}$. 
\end{lemma}
\begin{madef}[Modular extensions, modular diagrams] \label{defmoddiag}
Let $\L$ be a geometric lattice. 
\begin{itemize}
   \item  An \textit{extension of} $\L$ is the datum of a geometric lattice $\E$ together with an embedding $\imath: \L \rightarrow \E$ such that $\imath(\L)$ is an interval of the form $[\zero, F]$ with $F$ a flat of $\E$. In that case the element $F$ will be denoted $F_{\imath}$. 
   \item An extension $\imath:\L \hookrightarrow \E$ is called \textit{modular} if $F_{\imath}$ is modular in $\E$. 
    \item A \textit{modular diagram} $\Gamma$ of $\L$ is a triple $\Gamma = (\E, \imath, J)$ where $(\E, \imath)$ is a modular extension of $\L$ and $J$ is a linearly ordered multiset of atoms of $\E$. The multiset $J$ will be denoted $J_\Gamma$. 
     \item The degree of a modular diagram is defined by 
\begin{equation*}
    \deg \,(\E, \imath, J) = |J| - 2(\rk \medvee J - \rk (\medvee J)\wedge F).
\end{equation*}
\end{itemize}
\end{madef}
The definition of the modular extensions of a geometric lattice $\L$ may seem unnecessarily complicated, but this is the cleanest way of saying ``all the geometric lattices containing $\L$ as a modular flat'' in an equivariant-friendly way. By ``linearly ordered multiset of atoms of $\E$'' we mean a finite word in the alphabet $\At(\E).$ 
\begin{madef}\label{defMD}
Let $\L$ be a non-trivial geometric lattice. The $\Z$-graded vector space $\MD(\L)$ of modular diagrams of $\L$ is the vector space over $\Q$ freely spanned by the modular diagrams of $\L$ having degree as above, quotiented by the relations
    \begin{enumerate}[label={\arabic*)}]
        \item \label{rel1}$(\E, \imath,J) \sim -(\E, \imath, J')$ if $J$ and $J'$ differ by a transposition of their linear order. 
        \item \label{rel2}$(\E, \imath,J) \sim (\E', \imath', J')$ if there exists an embedding $\varphi: \E \hookrightarrow \E'$ such that $\rk(\E) = \rk(\E')$ and such that we have $\imath' = \varphi\circ \imath$ and $J' = \varphi(J)$. 
        \item \label{rel3}$(\E, \imath, J) \sim 0$ if $F_{\imath} \vee \bigvee J < \un_{\E}$.
        \item \label{rel4}$(\E, \imath, J) \sim 0$ if $\E$ is a product $\E_1\times\E_2$ with non-trivial $\E_1$ and $\E_2$, and we have $F_{\imath}\in \E_1\times \{\zero_{\E_2}\}$.
        \item \label{rel5} $(\E, \imath, J) \sim 0$ if there exists a modular element $F$ of $\E$ such that we have $F_{\imath} \leq F$ and such that the subset of elements of $J$ which are not below $F$ has cardinality $2$. 
    \end{enumerate}
\end{madef}
Note that we do not define the space of modular diagrams of a singleton, so we won't have access to $\Gra(1)$ for instance. This is because from the matroidal point of view we are unable to distinguish the case of one external vertex, that is, $\Gra(1)$ which is a silly object quasi-isomorphic to $\Q$, and the case of zero external vertex which we haven't defined and is quite complicated. Whenever we consider the space of modular diagrams of a geometric lattice it is implicitly assumed that this geometric lattice is non-trivial.

Note that relation \ref{rel1} implies that a modular diagram $\Gamma$ is zero in $\MD(\L)$ as soon as $J_{\Gamma}$ is not a set. A priori allowing $J_{\Gamma}$ to be a multiset will make our life easier later on. Note that in relation \ref{rel2} the embedding $\varphi$ could be an isomorphism. In the the braid case (geometric lattices associated to complete graphs) this will be responsible for the internal vertices (the vertices outside $F_{\imath}$) being unnumbered. Relation \ref{rel2} together with relation \ref{rel3} imply that we can always assume that we have $\At(\E) = F_{\imath}\cup J.$ Relations  \ref{rel3}, \ref{rel4} and \ref{rel5} are motivated by the results of the previous section. However, notice that relation \ref{rel4} is stronger than its geometric counterpart. This means that $\MD$ restricted to geometric lattices associated to complete graphs is not exactly $\Gra$. For instance the diagram 
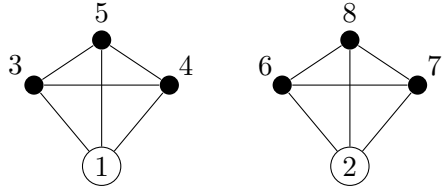
\begin{equation*}
    \begin{tikzpicture}[scale=0.60, baseline=14px]

    \node[circle, fill, black, inner xsep=0.9mm,inner ysep = 0.9mm ] (A) at (0.5,1.8) {};
    \node[circle, fill, black,inner xsep=0.9mm,inner ysep=0.9mm ] (B) at (3.5,1.8) {};
    \node[circle, fill, black,inner xsep=0.9mm,inner ysep=0.9mm ] (C) at (2,2.8) {};
    \node[circle, fill, black,inner xsep=0.9mm,inner ysep=0.9mm ] (D) at (6,1.8) {};
    \node[circle, fill, black,inner xsep=0.9mm,inner ysep=0.9mm ] (E) at (9,1.8) {};
    \node[circle, fill, black,inner xsep=0.9mm,inner ysep=0.9mm ] (F) at (7.5,2.8) {};
    \node[above left =-3pt of A, outer sep=0pt] {$3$};
    \node[above right=0pt of {(3.5,1.8)}, outer sep=0pt] {$4$};
    \node[above =3pt of {(2,2.8)}, outer sep=0pt] {$5$};
    \node[above left=0pt of {(6,1.8)}, outer sep=0pt] {$6$};
    \node[above right=0pt of {(9,1.8)}, outer sep=0pt] {$7$};
    \node[above =3pt of {(7.5,2.8)}, outer sep=0pt] {$8$};
    \node[circle, draw,inner xsep=0.7mm,inner ysep=0.7mm ] (G) at (2,0) {$1$};
    \node[circle, draw,inner xsep=0.7mm,inner ysep=0.7mm ] (H) at (7.5,0) {$2$};

    \draw[] (A) -- (G);
    \draw[] (B) -- (G);
    \draw[] (C) -- (G);
    \draw[] (A) -- (B);
    \draw[] (A) -- (C);
    \draw[] (B) -- (C);
    \draw[] (D) -- (H);
    \draw[] (E) -- (H);
    \draw[] (F) -- (H);
    \draw[] (D) -- (E);
    \draw[] (D) -- (F);
    \draw[] (E) -- (F);
\end{tikzpicture}
\end{equation*}
is not zero in $\Gra(2)$ but it is zero as a modular diagram of $\Pi_2$. To be more explicit, the modular extension we are considering is $\Pi_2 \hookrightarrow \Pi_8$ and the set $J$ is the set of atoms corresponding to the edges of the above graph (with any total order). The join of the atoms of $J$ is the flat $1345|2678$. By Relation \ref{rel2} we can restrict to $1345|2678 \cup 12|3|\ldots|8$, but that restriction is a direct product of $\Pi_2$ with the flat $1345|2678$ and so that modular diagram is zero in $\MD(\Pi_2)$ by relation \ref{rel4}. Let us describe another example in which $\Gra(n)$ and $\MD(\Pi_n)$ differ. Consider the following graphical diagram. 
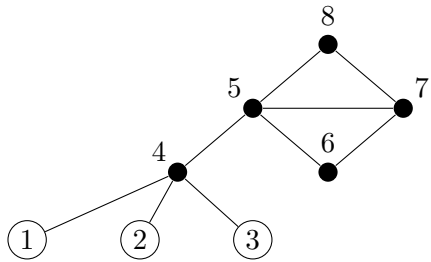
\begin{equation*}
    \begin{tikzpicture}[scale=0.50, baseline=14px]

    \node[circle, fill, black, inner xsep=0.9mm,inner ysep = 0.9mm ] (A) at (4,1.8) {};
    \node[circle, fill, black,inner xsep=0.9mm,inner ysep=0.9mm ] (B) at (6,3.5) {};
    \node[circle, fill, black,inner xsep=0.9mm,inner ysep=0.9mm ] (C) at (8,1.8) {};
    \node[circle, fill, black,inner xsep=0.9mm,inner ysep=0.9mm ] (D) at (10,3.5) {};
    \node[circle, fill, black,inner xsep=0.9mm,inner ysep=0.9mm ] (E) at (8,5.2) {};
    \node[above left =-3pt of A, outer sep=0pt] {$4$};
    \node[above left=-3pt of B, outer sep=0pt] {$5$};
    \node[above =0pt of C, outer sep=0pt] {$6$};
    \node[above right=-3pt of D, outer sep=0pt] {$7$};
    \node[above =0pt of E, outer sep=0pt] {$8$};
    \node[circle, draw,inner xsep=0.7mm,inner ysep=0.7mm ] (G) at (0,0) {$1$};
    \node[circle, draw,inner xsep=0.7mm,inner ysep=0.7mm ] (H) at (3,0) {$2$}; 
    \node[circle, draw,inner xsep=0.7mm,inner ysep=0.7mm ] (I) at (6,0) {$3$}; 

    \draw[] (I) -- (A);
    \draw[] (G) -- (A);
    \draw[] (H) -- (A);
    \draw[]  (A) -- (B);
    \draw[] (B) -- (C);
    \draw[] (C) -- (D);
    \draw[] (D) -- (E);
    \draw[] (E) -- (B);
    \draw[] (B) -- (D);
    \end{tikzpicture}
\end{equation*}
That picture defines a modular diagram of $\Pi_3$ with modular extension $\Pi_3 \hookrightarrow \Pi_8$ and set of atoms given by the edges of the graph (with any total order). By relation \ref{rel2} we can restrict to the atoms $1234|5678 \cup 45$ but that restriction is a direct product with one summand containing $\Pi_3$ and so by relation \ref{rel4} that modular diagram is zero in $\MD(\Pi_3)$. However the corresponding diagram in $\Gra(3)$ is not zero. In general, a diagram of $\Pi_n$ having a bridge between a subgraph containing $\Pi_n$ and another subgraph will be zero in $\MD(\Pi_n)$. In the rest of this article we will use the same notation $\Gamma = (\E, \imath, J)$ to denote a modular diagram and its equivalence class in $\MD(\L)$. If $\E$ is $\L$ and $\imath$ is the identity we will simply write $J$. If $\Gamma = (\E,\imath, J)$ is some modular diagram and $H$ is some element of $J$ we will denote $\Gamma\setminus H \coloneqq (\E, \imath, J \setminus H)$. Our first task is to define a differential on $\MD(\L).$ For this we will need the following classical result, known as the diamond isomorphism lemma. 
\begin{lemma}[\cite{birkhoff1940lattice}, Section IV.2]\label{lemmaisodiamond}
For any elements $F, F'$ in some geometric lattice $\E$, if $F$ is modular then taking the join with $F'$ defines an isomorphism of posets from $[F\wedge F', F]$ to $[F', F\vee F']$, with inverse given by taking the meet with $F$. 
\end{lemma}
We will denote this isomorphism by $\isomod_{F',F}$. We also have the following lemma. 
\begin{lemma}[{\cite[Lemma 2]{Stanley1971}\cite[Corollary 3.9]{Brylawski_1975}}]\label{lemmamodheredity}
For any elements $F,F'$ in a geometric lattice $\E$, if $F$ is modular in $\E$ then $F\vee F'$ is modular in $[F', \un]$ and $F\wedge F'$ is modular in $[\zero, F']$.
\end{lemma}
\begin{madeflemma}[Bridges, contractibility, contraction]
Let $\Gamma = (\E, \imath, J)$ be a modular diagram of $\L$. An atom $H \in J$ is called a \textit{bridge} if it is the only atom of $J$ not below some flat $F \geq F_{\imath}$. An atom $H \in J$ is called \textit{contractible} if $H$ is not a bridge and $H$ is not below $F_{\imath}$. The set of contractible atoms of $\Gamma$ will be denoted by $\Sc$. For any $H\in \Sc$, the \textit{contraction} $\Gamma / H$ of $\Gamma$ at $H$ is the modular diagram of $\L$ defined by $$\Gamma/H \coloneqq ([H, \un], H\vee \imath, H\vee (J\setminus \{H\}))$$ with linear order on $H\vee (J\setminus \{H\})$ induced by that on $J$. 
\end{madeflemma}
\begin{proof}
By Lemma \ref{lemmaisodiamond}, the morphism $H\vee \imath$ is an embedding with image $[H,H\vee F_{\imath}]\subset [H, \un]$ and by Lemma \ref{lemmamodheredity} the element $H\vee F_{\imath}$ is modular in $[H, \un]$, which implies that $([H, \un], H\vee \imath, H\vee J\setminus \{H\})$ is indeed a modular diagram of $\L$. 
\end{proof}
\begin{rmq}\label{rmkbridge}
Note that if a diagram $\Gamma = (\E, \imath, J)$ has a bridge $H$, unique atom of $J$ not below some flat $F \geq F_{\imath}$, then by relation \ref{rel3} in Definition \ref{defMD} we have $\Gamma \sim (\E_{|F \sqcup \{H\}}, \imath, J)$, but $\E_{|F \sqcup \{H\}}$ splits as the direct product $\E_{|F}\times [\zero, H]$ and so by relation \ref{rel3} we get $\Gamma \sim 0.$
\end{rmq}
\begin{ex}\label{excontraction}
Let us go back to the geometric lattice $\E$ described in Example \ref{exgeolatt}. As was noticed in Example \ref{exmodular} the inclusion $\imath: \L \coloneqq \E_{|\{a,b,c,d\}}\hookrightarrow \E$ is a modular extension. Let us consider the modular diagram $\Gamma = (\E, \imath, \{e\vartriangleleft b'\vartriangleleft c'\})$ of $\L$. We aim to compute the contraction $\Gamma/e$. Lemma \ref{lemmaisodiamond} tells us that we have an isomorphism $\vee_{e, \{a,b,c,d\}}:\L \simeq [e, \un]\subset \E$ given by taking the join with $e$. We have $e\vee b' = e\vee b$ which gives $\vee_{e, \{a,b,c,d\}}^{-1}(e\vee b')=b$, and $e\vee c' = e\vee c$ which gives $\vee_{e, \{a,b,c,d\}}^{-1}(e\vee c')=c$. This yields $$\Gamma /e = ([e,\un], \vee_{e, \{a,b,c,d\}}, \{e\vee b',e\vee c'\}) = (\L, \Id_{\L}, \{b,c\}).$$
Similar considerations give $\Gamma/b' = (\L, \Id_{\L}, \{b,a\})$ and $\Gamma/c' = (\L, \Id_{\L}, \{c,a\}).$
\end{ex}
\begin{rmq}\label{rmqdeg}
For all modular diagram $\Gamma$ and all contractible atom $H \in J_{\Gamma}^{contr}$ we have $\deg(\Gamma/H) = |J_{\Gamma}| - 1 - 2(\rk(\bigvee J) - \rk((\bigvee J) \wedge (F_{\imath} \vee H)))$, but by modularity of $F_{\imath}$ we have the equality $(\bigvee J)\wedge (F_{\imath} \vee H) = H \vee (\bigvee J \wedge F_{\imath})$ which gives $ \rk((\bigvee J)\wedge (F_{\imath} \vee H)) = \rk ((\bigvee J) \wedge F_{\imath}) + 1.$ This means that we have 
$$ \deg (\Gamma/H) =  |J_{\Gamma}| - 1 - 2(\rk(\bigvee J) - \rk((\bigvee J) \wedge F_{\imath}) - 1) = \deg(\Gamma) + 1.$$
\end{rmq}
We are now able to define the  differential of $\MD(\L)$.  
\begin{madeflemma}
The linear endomorphism $\d: \MD(\L)\rightarrow \MD(\L)$ taking value on modular diagrams $\Gamma = (\E, \imath, J= \{H_1 \vartriangleleft \cdots \vartriangleleft H_n\})$ 
\begin{equation*}
    \d(\Gamma) \coloneqq \sum_{H_i \in \Sc} (-1)^{i+1}\Gamma/H_i
\end{equation*}
is well-defined.
\end{madeflemma}
\begin{proof}
Let us check that the endomorphism $\d$ is compatible with relations \ref{rel1} to \ref{rel5} in Definition \ref{defMD}. For relation \ref{rel1} 1) one can see that if a modular diagram $\Gamma'$ can be obtained from another modular diagram $\Gamma$ by exchanging two atoms $H_p, H_{p+1}$ in the linear order of $J_{\Gamma}$ then for all $H \neq H_p, H_{p+1}$ we have $\Gamma'/H \sim -\Gamma/H$ because those two modular diagrams differ by a transposition of their linear order. Besides we obviously have $\Gamma/H_{p} = \Gamma'/H_{p}$ and $\Gamma/H_{p+1} = \Gamma'/H_{p+1}$ which finally gives $\d(\Gamma') \sim -\d(\Gamma)$. For relation \ref{rel2} the compatibility comes from the fact that if we have $\Gamma = (\E,\imath, J) \sim \Gamma' = (\E', \imath', J')$ via some embedding $\varphi:\E \hookrightarrow \E'$ then we have $\Gamma/H \sim \Gamma'/\varphi(H)$ via the embedding $\varphi_{|[H,\un]}.$ For relation \ref{rel3}, if for some modular diagram $\Gamma = (\E, \imath, J)$ we have $F_\imath \vee \bigvee J < \un_{\E}$ then one can see that for all $H \in \Sc$ we have $F_{H \vee \imath}\vee \bigvee(H\vee (J\setminus H))= F_{\imath}\vee \bigvee J < \un_{[H,\un]}$ which gives $\Gamma/H \sim 0$. For relation \ref{rel4}, if for some modular diagram $\Gamma = (\E, \imath, J)$ the geometric lattice $\E$ is a non-trivial product $\E'\times \E''$ with $\E'$ containing $F_{\imath}$, then for any $H\in\Sc$ the interval $[H, \E]$ is either isomorphic to $[H, \un_{\E'}]\times\E''$ or to $\E'\times [H, \un_{\E''}]$, which are both non-trivial products. Finally, for relation \ref{rel5} if a modular diagram $\Gamma = (\E, \imath, J)$ is such that we have exactly two atoms $H,H' \in J$ not below some modular flat $F\geq F_{\imath},$ then by relation \ref{rel3} we can assume that we have $F\vee H_1\vee H_2 = \un.$ First consider the case $\rk (F\vee H_1 \vee H_2) - \rk(F) = 1.$ In that case for all $H\in \Sc$ different from $H_1$ and $H_2$, the flats $H\vee H_1$ and $H\vee H_2$ are the only atoms of $[H,\un]$ not below $F$ and so by Remark \ref{rmkbridge} we have $\Gamma/H \sim 0$. Besides, by the diamond isomorphism lemma \ref{lemmaisodiamond} applied to $F$ the contractions $\Gamma/H_1$ and $\Gamma/H_2$ cancel out in $\d(\Gamma).$ In the case $\rk (F\vee H_1 \vee H_2) - \rk(F) = 2$ the contractions $\Gamma/H$ are equivalent to $0$ for all $H \neq H_1,H_2$ for the same reasons as above. Besides, the contraction $\Gamma/H_1$ is zero because $H_1\vee H_2$ is the unique atom of $H_1 \vee (J\setminus H_1)$ which is not below the modular flat $F\vee H_1 \in [H_1, \un]$, and the same goes for $\Gamma/H_2$. 
\end{proof}
\begin{ex}
Going back to the modular diagram $\Gamma$ considered in Example \ref{excontraction}, the computations in that example give 
\begin{equation*}
    \d(\Gamma) = (\L, \Id_{\L}, \{b,c\}) - (\L, \Id_{\L}, \{b,a\}) + (\L, \Id_{\L}, \{c,a\}).
\end{equation*}

\end{ex}
By Remark \ref{rmqdeg} the endomorphism $\d$ has degree $+1$. Besides, we have the following lemma. 
\begin{lemma}
The endomorphism $\d$ satisfies the equation $\d^2 = 0$.
\end{lemma}
\begin{proof}
Let $\Gamma = (\E, \imath, J= \{H_1 \vartriangleleft \cdots \vartriangleleft H_n\})$ be a modular diagram of $\L$, and let $H_i \vartriangleleft H_j \in \Sc$ be two contractible atoms. If the rank $2$ flat $H_i\vee H_j$ contains strictly more than two elements in $J$, then the term $(-1)^i(-1)^{j-1}(\Gamma/{H_i})/(H_i\vee H_j)$ appearing in the expansion of $\d^2(\Gamma)$ is $0$ because $J_{\Gamma/{H_i}}$ is not a set. Otherwise, that term is cancelled by $(-1)^j(-1)^i(\Gamma/H_j)/H_i\vee H_j$.
\end{proof}
Our next task is to define a commutative product on $\MD(\L)$ which is compatible with our newly defined differential. As we have seen in the previous section, at the level of modular extensions of $\L$ this product should be given by ``gluing extensions along $\L$''. In general gluing matroids along a minor, usually called 'amalgamation', is a non-trivial process. In particular, the gluing may not exist, or there may be several non-canonical gluings (see \cite[Section 11.4]{oxley}). However, in our case everything works out nicely because we are gluing along a modular flat. In this particular situation the gluing is called a \textit{generalized parallel connection} (see \cite[Section 5]{Brylawski_1975}). In this article we will prefer the more suggestive terminology of \textit{pushout}. This construction is indeed a pushout in a well-chosen category of geometric lattices but this will not play a role in this article. 
\begin{madef}[Pushout]\label{defpushout}
Let $\imath_1: \L \hookrightarrow \E_1$ and $\imath_2: \L \hookrightarrow \E_2$ be two modular extensions of $\L$. The \textit{pushout of $\E_1$ and $\E_2$ along $\L$}, denoted by $\E_1\cup_{\L}\E_2$, is the subposet of $\E_1\times \E_2$ consisting of elements $(F_1,F_2)\in \E_1\times \E_2$ such that we have the equality $\imath_{1}^{-1}(F_1\wedge F_{\imath_1}) = \imath_{2}^{-1}(F_2\wedge F_{\imath_2}).$ We denote by $\imath_{12}: \L \hookrightarrow \E_1\cup_{\L}\E_2$ the embedding sending a flat $F\in \L$ to $(\imath_{1}(F), \imath_{2}(F))$.
\end{madef}
\begin{ex}
If $\L = \{\star\}$ we simply get $\E_1\cup_{\L}\E_2 = \E_1\times \E_2,$ which is indeed the coproduct in the category of finite lattices for instance.
\end{ex}
\begin{ex}
Consider the geometric lattice $\L$ described in Example \ref{exgeolattres}. That geometric lattice is (isomorphic to) the pushout of the restrictions $\L_{|\{a,b,c,d\}}$ and $\L_{|\{a,b',c',d'\}}$ along the flat $\{a\}$ (which is modular in those restrictions because an atom is always modular). 
\end{ex}
\begin{ex}
If $\E_1$ and $\E_2$ are the geometric lattices associated to some graphs $G_1$, $G_2$, and $\L$ is the geometric lattice associated to some common modular subgraph $G$ then the pushout $\E_1\cup_{\L}\E_2$ is the geometric lattice associated to the graph obtained by gluing $G_1$ and $G_2$ along $G$.
\end{ex}
\begin{rmq}\label{rmkatomspushout}
    The atoms of a pushout $\E_1\cup_{\L}\E_2$ are the elements of the form $(H,\zero)$ for $H\in \At(\E_1)\setminus \iota_1(\L)$, the elements of the form $(\zero, H)$ for $H \in \At(\E_2)\setminus \iota_2(\L)$, and the elements of the form $(\iota_1(H), \iota_2(H))$ for $H \in \At(\L).$ We can thus naturally identify $\At(\E_1\cup_{\L}\E_2)$ with $\At(\E_1)\cup_{\At(\L)}\At(\E_2).$
\end{rmq}
\begin{rmq}\label{rmkembedpushout}
For any pushout $\E_1\cup_{\L}\E_2$ along $\imath_1, \imath_2$, we have an embedding $\imath_{\E_1}: \E_1\hookrightarrow \E_1\cup_{\L}\E_2$ defined by sending a flat $F\in \E_1$ to $(F, \imath_2\imath_1^{-1}(F\wedge F_{\imath_1})) \in \E_1\cup_{\L}\E_2$, and similarly for $\E_2$ by symmetry. Those embeddings give isomorphisms $\E_1 \simeq (\E_1\cup_{\L}\E_2)_{|\At(\E_1)}$ and $\E_2 \simeq (\E_1\cup_{\L}\E_2)_{|\At(\E_2)}$
\end{rmq}
We summarize the relevant properties of pushouts in the following lemma. We refer to \cite[Proposition 5.10]{Brylawski_1975} and \cite[Proposition 11.4.14]{oxley} for the proofs. 
\begin{lemma}\label{lemmapushout}
Let $\imath_1:\L \hookrightarrow \E_1$ and $\imath_2: \L \hookrightarrow \E_2$ be two modular extensions.
\begin{enumerate}[label={\arabic*)}]
    \item The poset $\E_1\cup_{\L} \E_2$ is a geometric lattice with rank given by $$\rk((F_1, F_2)) = \rk_{\E_1}(F_1) + \rk_{\E_2}(F_2) - \rk_{\L}(\imath_1^{-1}(F_1\wedge F_{\imath_1})).$$ Moreover, this is still true even if only one of $\imath_1$ and $\imath_2$ is modular. 
    \item The embedding $\imath_{\E_1}: \E_1 \hookrightarrow \E_1\cup_{\L}\E_2$ defined in Remark \ref{rmkembedpushout} is a modular extension. Moreover this is true even if only $\imath_2$ is modular.  
    \item The embedding $\imath_{12}$ is a modular extension.
    \item The modular extensions $\E_1\cup_\L\E_2$ and $\E_2\cup_{\L}\E_1$ are isomorphic. 
    \item If we have a third modular extension $\imath_3: \L \hookrightarrow \E_3$ then the modular extensions $(\E_1\cup_{\L} \E_2)\cup_{\L}\E_3$ and $\E_1\cup_\L(\E_2\cup_{\L}\E_3)$ are isomorphic. 
    \item If $H$ is an atom of $\E_1$ which is not below $F_{\imath_1}$ then the modular extension $(\E_1\cup_{\L}\E_2)_{\geq H}$ (with embedding $((H,\zero)\vee (\imath_1\times \imath_2))$) and the modular extension $(\E_1)_{\geq H}\cup_{\L}\E_2$ (with embedding $\L\hookrightarrow (\E_1)_{\geq H}$ given by $H\vee\imath_1$) are isomorphic.  
    \item For any $(F_1, F_2)\in \E_1\cup_\L\E_2$ the interval $[\zero, (F_1, F_2)]$ is isomorphic to $[\zero, F_1]\cup_{[\zero, \imath_{1}^{-1}(F_1\wedge F_{\imath_1})]}[\zero, F_2]$, and the interval $[(F_1, F_2), \un]$ is isomorphic to $[F_1, \un]\cup_{[\imath_{1}^{-1}(F_1\wedge F_{\imath_1}), \un]}[F_2,\un].$ 
    \item For any subset $S \subset \At(\E_1\cup_{\L}\E_2)$ containing $\At(\L)$, the restriction $(\E_1\cup_{\L}\E_2)_{|S}$ is isomorphic to the pushout $(\E_1)_{|S\cap \At(\E_1)}\cup_{\L}(\E_2)_{|S\cap \At(\E_2)}$.
\end{enumerate}
\end{lemma}
We will also need the following two lemmas.
\begin{lemma}[{\cite[Proposition 5.9]{Brylawski_1975}}]\label{lemmadetectionpushout}
Let $F_1$ and $F_2$ be two flats both containing a third flat $F$ in a geometric lattice $\E$. If $F$ is modular in $[\zero, F_1]$ then $\E$ is isomorphic to the pushout $[\zero, F_1]\cup_{[\zero,F]} [\zero, F_2]$ if and only if we have $[F,\un] \simeq [F,F_1]\times [F,F_2]$.   
\end{lemma}
\begin{lemma}\label{lemmamodularpushout}
Let $\imath_1:\L\hookrightarrow \E_1$ and $\imath_2:\L \hookrightarrow \E_2$ be two modular extensions of $\L$ and let $F$ be a flat of $\E_1$ greater than or equal to $F_{\imath_1}$. If $F$ is modular in $\E_1$ then $(F, \un_{\E_2})$ is modular in $\E_1\cup_{\L}\E_2$.
\end{lemma}
\begin{proof}
Let $(F_1,F_2)$ be a flat of $\E_1\cup_L\E_2$. We have 
\begin{align*}
    \rk((F_1,F_2)\vee (F, \un_{\E_2})) + \rk((F_1,F_2)\wedge (F, \un_{\E_2})) &= \rk((F_1\vee F, \un_{\E_2})) + \rk((F_1\wedge F, F_2)) \\
    &= \rk (F_1\vee F) + \rk(\E_2) - \rk(F_{\imath_1})\\
    & \,\,\,\,\,\,\,+ \rk(F_1\wedge F) + \rk(F_2) - \rk(F_1\wedge F\wedge F_2) \\
    &= \rk(F_1) + \rk(F) + \rk(\E_2) - \rk(F_{\imath_1})\\
    & \,\,\,\,\,\,\, + \rk(F_2) - \rk(F_1\wedge F\wedge F_2) \\
    &= \rk((F, \un_{\E_2})) + \rk((F_1,F_2)),
\end{align*}
where the third equality comes from the modularity of $F$ in $\E_1$. 
\end{proof}
We now proceed to define the commutative product on $\MD(\L)$. If $J_1$ and $J_2$ are two linearly ordered multisets in some common set, we denote by $J_1\cup J_2$ the multiset union of $J_1$ and $J_2$ with linear order given by concatenating the orders on $J_1$ and $J_2$, putting the elements of $J_1$ before the elements of $J_2$. 
\begin{madeflemma}
The bilinear endomorphism $\bullet:\MD(\L)\otimes \MD(\L) \rightarrow \MD(\L)$ taking value on modular diagrams 
\begin{equation*}
    (\E_1, \imath_1, J_1)\bullet(\E_2, \imath_1, J_2) \coloneqq (\E_1 \cup_{\L}\E_2, \imath_{12}, J_1\cup J_2),
\end{equation*}
is well-defined. 
\end{madeflemma}
Here we are taking the union $J_1\cup J_2$ in $\At(\E_1\cup_{\L}\E_2) = \At(\E_1)\cup_{\At(\L)}\At(\E_2)$ (see Remark \ref{rmkatomspushout}).
\begin{proof}
We must check that the bilinear endomorphism $\bullet$ is compatible with relations \ref{rel1} to \ref{rel5} in Definition \ref{defMD}. For relation \ref{rel1} it is obvious that if two modular diagrams $\Gamma_1,\Gamma_1'$ of $\E_1$ differ by a transposition of their linear order then $\Gamma_1\bullet \Gamma_2,\Gamma'_1 \bullet \Gamma_2$ will do too for any modular diagram $\Gamma_2$ of $\E_2$, which gives the result by symmetry. For relation \ref{rel2}, if $\Gamma_1, \Gamma_1'$ are two modular diagrams which are equivalent via some embedding $\varphi$, then by Lemma \ref{lemmapushout} item 8) the modular diagrams $\Gamma_1\bullet \Gamma_2$ and $\Gamma_1'\bullet \Gamma_2$ are also equivalent via some embedding for all $\Gamma_2$, which gives the result by symmetry. For relation \ref{rel3} if some modular diagram $\Gamma_1$ of $\E_1$ is such that we have $F_{\imath_1}\vee J_{\Gamma_1} < \un_{\E_1}$ then for all modular diagram $\Gamma_2$ in $\E_2$ we will have $(F_{\imath_1}, F_{\imath_2})\vee (J_{\Gamma_1}\cup J_{\Gamma_2}) \leq (F_{\imath_1}\vee J_{\Gamma_1}, \un_{\E_2}) < \un_{\E_1\cup_{\L}\E_2},$ which implies $\Gamma_1\bullet \Gamma_2 \sim 0$, which gives the result by symmetry. For relation \ref{rel4} assume we have two modular diagrams $\Gamma_1= (\E_1,\imath_1, J_1), \Gamma_2 = (\E_2, \imath_2, J_2)$ such that $\E_1$ is some non-trivial product $\E'_1\times \E''_1$ with $\E'_1$ containing $F_{\imath_1}$. One can check that the pushout $\E_1\cup_{\L}\E_2$ is isomorphic to the product $(\E'_1\cup_{\L}\E_2)\times \E''_1$ which implies $\Gamma_1\bullet \Gamma_2 \sim 0$. For relation \ref{rel5} if a modular diagram $\Gamma_1$ of $\E_1$ is such that $J_{\Gamma_1}$ contains exactly two atoms $H_1, H_2$ not below some modular flat $F \geq F_{\Gamma_1}$ then for all modular diagram $\Gamma_2$ in $\E_2$, by Lemma \ref{lemmamodularpushout} the element $(F_{\Gamma_1}, \un_{\E_2})$ is a modular flat of $\E_1\cup_\L\E_2$ greater than or equal to $(F_{\Gamma_1}, F_{\Gamma_2})$, and $H_1,H_2$ viewed as atoms of $\E_1\cup_{\L}\E_2$ are the unique atoms of $J_1\cup J_2$ not below that modular flat, and so we get $\Gamma_1\bullet \Gamma_2 \sim 0.$ 
\end{proof}
\begin{lemma}
The binary product $\bullet$ is graded commutative, associative and unitary. 
\end{lemma}
\begin{proof}
Commutativity is a direct consequence of Lemma \ref{lemmapushout} item 4), and associativity is a direct consequence of Lemma \ref{lemmapushout} item 5). The unit is given by the empty modular diagram with modular extension $\L \xrightarrow{\Id} \L$. 
\end{proof}
Next we must prove that the product and the differential on $\MD(\L)$ are compatible. 
\begin{lemma}\label{lemmaleibniz}
The binary product $\bullet$ satisfies the Leibniz identity
\begin{equation*}
    \d(\Gamma_1\bullet \Gamma_2) = \d(\Gamma_1)\bullet \Gamma_2 + (-1)^{\deg(\Gamma_1)}\Gamma_1\bullet \d(\Gamma_2). 
\end{equation*}
\end{lemma}
\begin{proof}
We can identify $J_{\Gamma_1\bullet \Gamma_2}^{contr} = J_{\Gamma_1}^{contr}\sqcup J_{\Gamma_2}^{contr}$. Furthermore by Lemma \ref{lemmapushout} item 6), for all $H$ in $J_{\Gamma_1}^{contr}$ (resp. in $J_{\Gamma_2}^{contr}$) we have $(\Gamma_1\bullet \Gamma_2)/H = (\Gamma_1/H)\bullet \Gamma_2$ (resp. $(\Gamma_1\bullet \Gamma_2)/H = \Gamma_1\bullet (\Gamma_2/H))$. This gives   
\begin{align*}
    \d(\Gamma_1\bullet \Gamma_2) &= \sum_{H_i \in J_{\Gamma_1\bullet \Gamma_2}^{contr}}(-1)^{i+1}(\Gamma_1\bullet \Gamma_2)/H_i \\
    &= \sum_{H_i \in J_{\Gamma_1}^{contr}}(-1)^{i+1}(\Gamma_1/H_i)\bullet \Gamma_2 + \sum_{H_j \in J_{\Gamma_2}^{contr}}(-1)^{|J_{\Gamma_1}|+j+1}\Gamma_1\bullet (\Gamma_2/H_j) \\
    &= \d(\Gamma_1)\bullet \Gamma_2 + (-1)^{\deg(\Gamma_1)}\Gamma_1\bullet \d(\Gamma_2).
\end{align*}
\end{proof}
From now on the notation $\MD(\L)$ refers to the cdga of modular diagrams with differential $\d$ and commutative product $\bullet$. The next proposition shows that $\MD(\L)$ is graded by the elements of $\L$, as is often the case for algebraic invariants of matroids. 
\begin{madefprop}\label{defgrading}
Let $\L$ be a geometric lattice. The $\L$-grading of a modular diagram $(\E, \imath, J) \in \MD(\L)$ is the flat of $\L$ defined by $\gr_{\L}(\E, \imath, J) \coloneqq \imath^{-1}((\bigvee J) \wedge F_{\imath}) \in \L$. If we denote by $\MD(\L, F)$ the subspace of modular diagrams of grading $F \in \L$, we have a decomposition of dg-complexes $\MD(\L) = \bigoplus_{F \in \L} \MD(\L, F)$. This decomposition is compatible with the commutative product $\bullet$ in the sense that we have $\MD(\L, F_1)\bullet \MD(\L,F_2) \subset \MD(\L, F_1\vee F_2).$ Furthermore for all $F\in \L$ we have an isomorphism of dg-complexes $\MD(\L, F) \simeq \MD([\zero, F], \un).$
\end{madefprop}
\begin{proof}
For all modular diagram $\Gamma = (\E, \imath, J)$ and for all contractible atom $H \in \Sc$, by modularity of $F_{\imath}$ we have 
\begin{align*}
    \gr_{\L}(\Gamma/H) & = (H\vee \imath)^{-1}((\medvee J) \wedge (H \vee F_{\imath})) \\
    &= (H\vee \imath)^{-1}(H \vee ((\medvee J) \wedge F_{\imath})) \\
    &= \imath^{-1}((\medvee J) \wedge F_{\imath}) \\
    &= \gr_{\L}(\Gamma).
\end{align*}
This implies that for all flat $F\in \L$ the subspace $\MD(\L,F)$ is a sub dg-complex of $\MD(\L)$, which gives the claimed decomposition. For all modular diagrams $\Gamma_1 = (\E_1, \imath_1, J_1), \Gamma_1 = (\E_1, \imath_1, J_1)$ of $\L$ we have
\begin{align*}
    F_{\imath_{12}}\wedge (\medvee (J_1 \cup J_2)) &= F_{\imath_{12}}\wedge (\medvee J_1 \vee \medvee J_2) \\
    &= F_{\imath_{12}}\wedge F_1\wedge (\medvee J_1 \vee \medvee J_2) \\
    &= F_{\imath_{12}}\wedge ( \medvee J_1 \vee (F_1\wedge \medvee J_2)) \\ 
    &= (F_{\imath_{12}} \wedge \medvee J_1)\vee (F_{\imath_{12}}\wedge \medvee J_2).
\end{align*}
Applying $\imath_{12}^{-1}$ we get $\gr_{\L}(\Gamma_1 \bullet \Gamma_2) = \gr_{\L}(\Gamma_1) \vee  \gr_{\L}(\Gamma_1),$ which gives the claimed inclusion. For all $F \in \L$ and all modular diagram $(\E, \imath, J)$ of grading $F$, by Lemma \ref{lemmamodheredity} the triple $(\E_{|\iota(F)\cup J}, \iota_{|F}, J)$ is a modular diagram of $[\zero,F]$, which has grading $\un_{[\zero,F]}$. This defines a morphism $\MD(\zero, F) \rightarrow \MD([\zero, F], \un).$ This morphism has an inverse given by sending a modular diagram $(\E, \imath, J)\in \MD([\zero, F], \un)$ to the modular diagram $(\L\cup_{[\zero, F]} \E, \imath_{\L}, J)\in \MD(\L, F)$. Note that even though $F$ is a priori not modular in $\L$, by assumption $F_{\imath}$ is modular in $\E$ which is enough to define the pushout $\L\cup_{[\zero, F]}\E$ (see Lemma \ref{lemmapushout} item 1)). Furthermore, the embedding $\L \hookrightarrow \L\cup_{[\zero, F]} \E$ is a modular extension by Lemma \ref{lemmapushout} item 2). 
\end{proof}
\begin{rmq}\label{rmkgradingzero}
Let $\Gamma = (\E, \imath, J) \in \MD(\L)$ be a modular diagram of $\L$-grading $\zero$ and such that $J$ is non-empty. If we denote $F = \bigvee J$, by relation \ref{rel2} we can assume $\At(\E) = \At(\L)\sqcup F$. By modularity of $F_{\imath}$ we have $\rk (\E) = \rk(F_{\imath}) + \rk(F)$, which implies that $\E$ decomposes as a direct product $\E \simeq \L\times [\zero, F]$. By relation \ref{rel4} this implies that we have $\Gamma \sim 0$. In other words, we have just proved that $\MD(\L, \zero)$ is isomorphic to $\Q$ concentrated in degree $0$, generated by the empty modular diagram. 
\end{rmq}
\begin{prop}
Let $\L_1$ and $\L_2$ be two geometric lattices. We have an isomorphism of cdgas $\MD(\L_1)\otimes \MD(\L_2) \simeq \MD(\L_1\times \L_2)$.
\end{prop}
\begin{proof}
One can easily check that if $\imath_1:\L_1\hookrightarrow\E_1$ and $\imath_2:\L_2 \hookrightarrow \E_2$ are two modular extensions, then the product embedding $\imath_1\times \imath_2:\L_1\times \L_2 \hookrightarrow \E_1 \times \E_2$ is also a modular extension. This allows us to define a morphism of cdgas $\MD(\L_1)\otimes \MD(\L_2) \rightarrow \MD(\L_1\times \L_2)$ by sending a pair of modular diagrams $(\E_1, \imath_1, J_1), (\E_2, \imath_2, J_2)$ to the modular diagram $(\E_1\times \E_2, \imath_1\times \imath_2, J_1 \sqcup J_2)$. By \cite[Corollary 3.16]{Brylawski_1975} any modular extension of $\L_1\times \L_2$ is a product of modular extensions of $\L_1$ and $\L_2$ respectively, which implies that the above morphism is an isomorphism. 
\end{proof}
\begin{prop}\label{propfree}
For any geometric lattice $\L$, the cdga $\MD(\L)$ is free as a graded commutative algebra. 
\end{prop}
\begin{proof}
Let us say that a modular extension $\imath: \L\hookrightarrow \E$ is irreducible if is not a pushout of two non-trivial modular extensions of $\L$. By Lemma \ref{lemmadetectionpushout} and Proposition \ref{propdecir} applied to $[F_{\imath}, \un]$ any modular extension decomposes uniquely as a pushout of irreducible modular extensions. This implies that the commutative graded algebra $\MD(\L)$ is free with generators given by equivalence classes of modular diagrams $\Gamma = (\E,\imath, J)$ such that $\imath$ is irreducible.
\end{proof}
We now turn to proving that in some good cases, the cdga $\MD(\L)$ is a model of $\OS(\L)$, meaning that we have a quasi-isomorphism of cdgas $\MD(\L) \simeq \OS(\L)$. The following proposition shows that we always have a well-defined morphism between those two cdgas.  
\begin{prop}\label{propmor}
For any geometric lattice $\L$ we have a morphism of cdgas $$\I_{\L}:\MD(\L) \rightarrow (\OS(\L),\d = 0)$$  defined on modular diagrams by $\I_{\L}(\E, \imath, J) \coloneqq e_{\imath^{-1}(J)}$ if no element of $J$ is contractible, and $\I_{\L}(\E, \imath, J) \coloneqq 0$ otherwise. 
\end{prop}
If the geometric lattice $\L$ is clear from the context we will simply write $\I$.
\begin{proof}
The fact that $\I$ is compatible with relations \ref{rel1} to \ref{rel5} in Definition \ref{defMD} is immediate, and so is the fact that $\I$ is a morphism of algebras. For the compatibility with the differential we need to show that for all modular diagrams $\Gamma \in \MD(\L)$ we have the equality $\I(\d(\Gamma)) = 0$. Let us denote $\Gamma = (\E, \imath, J)$. If no element of $J$ is contractible then we have $\d(\Gamma) = 0$ so the result follows. Furthermore, if $\rk \medvee J - \rk (\medvee J) \wedge F_{\imath}$ is strictly greater than $1$, then $\d(\Gamma)$ is a sum of modular diagrams all containing at least one contractible atom, which is sent to zero by $\I$. Let us assume that we have $\rk \medvee J - \rk (\medvee J) \wedge F_{\imath} = 1$. Notice that up to a sign we have the equality $$ \Gamma = \pm (\L, \Id_{\L}, J\setminus \Sc)\bullet (\E, \imath, \Sc)$$ and so by the Leibniz identity (Lemma \ref{lemmaleibniz}) we can assume $J = \Sc$. We now prove the result by induction on the cardinality of $J$. To simplify notations let us assume that $\imath$ is an inclusion. If $J$ has cardinality $3$ (the base case) and $\rk \bigvee J = 2$ then $\d(\Gamma)$ is already $0$. Otherwise if $\rk \bigvee J = 3$, let us denote $J = \{H_1,H_2,H_3\}$, $F = (\medvee J) \wedge F_{\imath}$, and $H_{ij} = (\vee_{H_i, F_{\imath}})^{-1}(H_i\vee H_j)$ for all $i\neq j \leq 3$. Note that for all $i \neq j \leq 3$ we have $H_{ij} = H_{ji} = (H_i \vee H_j)\wedge F.$ We can then compute $\I(\d(\Gamma)) = \pm\delta(e_{H_{12}}e_{ H_{13}}e_{H_{23}})$ which is zero since $\{H_{12}, H_{13}, H_{23}\}$ is a multiset of cardinality $3$ in the rank $2$ flat $F$. Next, assume that $J$ has cardinality $n$ for some integer $n$. Denote by $H_1 \vartriangleleft \cdots \vartriangleleft H_n$ the elements of $J$, which we can assume to be distinct. Denote by $H$ the element $(H_1\vee H_2)\wedge F_{\imath}$, which is an atom of $\L$. One can compute
\begin{align*}
    \d(\Gamma) - H\bullet \d(\Gamma\setminus H_2) + & H\bullet \d(\Gamma\setminus H_1) = \Gamma/H_1 - \Gamma/H_2 + \sum_{3 \leq k \leq n}(-1)^{k+1}\Gamma/H_k \\  
    &- H\bullet ((\Gamma\setminus H_2)/H_1) + \sum_{3\leq k\leq n}(-1)^{k+1} H\bullet ((\Gamma\setminus H_2)/H_k)\\
    & + H\bullet ((\Gamma\setminus H_1)/H_1) + \sum_{3 \leq k \leq n}(-1)^k H\bullet ((\Gamma\setminus H_1)/H_k) \\ 
    & = \sum_{3 \leq k \leq n }(-1)^{k+1}((H_1\vee H_k)\wedge F_{\imath})\bullet ((H_2\vee H_k)\wedge F_{\imath})\bullet((\Gamma\setminus H_1\setminus H_2)/H_k) \\
    & + \sum_{3 \leq k \leq n }(-1)^{k+1}H\bullet ((H_1\vee H_k)\wedge F_{\imath})\bullet((\Gamma\setminus H_1\setminus H_2)/H_k)\\ 
    &- \sum_{3 \leq k \leq n }(-1)^{k+1}H\bullet ((H_2\vee H_k)\wedge F_{\imath})\bullet((\Gamma\setminus H_1\setminus H_2)/H_k) \\
    &= \sum_{3 \leq k \leq n }(-1)^{k+1}\d(\{H_1, H_2, H_k\}) \bullet((\Gamma\setminus H_1\setminus H_2)/H_k),
\end{align*}
which gives the induction step. 
\end{proof}
Note that the above proof shows that $\I$ factors through the quadratic envelope of $\OS(\L)$ and so as soon as $\OS(\L)$ is not quadratic, the morphism $\I$ cannot be a quasi-isomorphism. For instance, if we consider the graphical geometric lattice $\L_{C_4}$ associated to the $4$-cycle $C_4$, with edges labeled from $1$ to $4$, then we have 
$$\OS(\L_{C_4}) = \frac{\bigwedge^{\bullet}[e_1, e_2, e_3, e_4]}{(e_1e_2e_3 - e_1e_2e_4 + e_1e_3e_4 - e_2e_3e_4)},$$ in which case $\I_{\L_{C_4}}$ is not a quasi-isomorphism. Determining which Orlik--Solomon algebras are quadratic is an interesting open question (see \cite[Section 4]{Falk_2001} \cite[Section 6]{Yuzvinsky_2001} for instance) which we will avoid in this article. Instead we will focus on a particular class of geometric lattices which are known to have quadratic (and even Koszul) Orlik--Solomon algebras, called supersolvable geometric lattices. To motivate their definition let us go back to braid arrangements. Recall from Section \ref{secgeo} that complements of braid arrangements admit a sequence of fibrations
\begin{equation*}
    \A_{\Br_n} \twoheadrightarrow \cdots \twoheadrightarrow \A_{\Br_2}\simeq \Cx^* \twoheadrightarrow \{*\}
\end{equation*}
given by forgetting points one after the other. This has strong consequences on the homotopy type of $\A_{\Br_n}$, for instance this implies that $\A_{\Br_n}$ is a $K(\pi, 1)$ space for all $n$. As explained in Section \ref{secgeo} the above fibrations can be seen at the combinatorial level, by noticing that the partitions 
\begin{equation*}
    12|3|\ldots|n < 123|4|\ldots|n < \cdots < 12\ldots n-1|n \in \Pi_n 
\end{equation*}
are modular elements of $\Pi_n$. This situation was generalized by Stanley \cite{Stanley_1972} as follow.
\begin{madef}[Supersolvable lattice]
    A geometric lattice is called \textit{supersolvable} if it admits a maximal chain of modular elements. 
\end{madef}
The terminology comes from supersolvable groups, whose lattice of subgroups admit a maximal chain of normal (and hence modular) subgroups. 
\begin{ex}
The geometric lattice $\L_G$ associated to a graph $G$ is supersolvable if and only if $G$ is chordal, meaning that every circuit in $G$ of length greater than or equal to $4$ has a chord. For instance the $4$-cycle $C_4$ does not have chords and so $\L_{C_4}$ is not supersolvable.
\end{ex}
The rest of this section will be devoted to proving the following theorem. 
\begin{theo}\label{theoqiso}
If a geometric lattice $\L$ is supersolvable then $\I_{\L}$ is a quasi-isomorphism. 
\end{theo}
One could vaguely sum up this result as saying that having enough modular flats implies having enough modular extensions. We must set the stage before starting the proof. Currently, our main issue is that we do not have any general method for producing modular extensions. Going back to partition lattices, in that case modular extensions of $\Pi_n$ are simply provided by embedding $\Pi_n$ in a bigger partition lattice $\Pi_n \hookrightarrow \Pi_{n+1}$. In terms of graphs, one can produce the $n+1$-th complete graph $K_{n+1}$ out of $K_n$ by first gluing two copies of $K_n$ along $K_{n-1}$ and then adding an edge between the two apexes. In matroidal terms we already know how to glue matroids (see Definition \ref{defpushout}) and there exists an explicit method for adding an element to a matroid in a controlled way, which we briefly describe now. We refer to \cite[Section 7.2]{oxley} for more details. The general idea is to freely add the element and then perform a so-called truncation along a modular cut. 
\begin{madef}[Modular cut]
Let $\L$ be a geometric lattice. A subset $\M \subset \L$ is called a \textit{modular cut} if it is upward-closed and if for all $F_1, F_2\in \M$ such that $\rk(F_1) + \rk(F_2) = \rk(F_1\wedge F_2) + \rk (F_1\vee F_2)$ we have $F_1 \wedge F_2 \in \M$. 
\end{madef}
\begin{madef}[Truncation]
Let $\M$ be a modular cut of a geometric lattice $\L$. The \textit{truncation of $\L$ along $\M$} is the geometric lattice 
$$ \Tr_{\M}(\L) \coloneqq \L \setminus \{F \in \L\setminus\M \, | \, F \textrm{ is covered by an element of } \M \}.$$
\end{madef}
One can easily check that if $\M$ is a modular cut of a geometric lattice $\L$ then the subset
 $\M\times \{e\} \subset \L\times \{\zero < e\}$ (for some symbol $e$) is a modular cut of $\L\times \{\zero < e\}$. 
\begin{madef}[Single-element extension]
Let $\M$ be a modular cut of a geometric lattice $\L$, and $e$ some symbol. The \textit{single-element extension of $\L$ along $\M$} is the geometric lattice $$\L\cup_{\M} e \coloneqq \Tr_{\M\times e}(\L\times \{\zero < e\}).$$
\end{madef}
Unless $\M$ contains an atom (which is a degenerate case), the atoms of $\L\cup_{\M} e$ are the atoms of $\L$ together with the symbol $e$. The map $\imath: \L \rightarrow \L\cup_{\M} e$ sending a flat $F$ to $(F, \zero)$ if $F\notin \M$ and to $(F, e)$ otherwise, is an embedding. Notice that by definition a flat $\imath(F)$ contains $e$ if and only if $F$ is in $\M$. In other words the choice of $\M$ dictates in which flats of $\L$ we add $e$. Going back to our example of complete graphs, the edge between the two apexes is added to every flat of $K_n\cup_{K_{n-1}} K_n$ containing two copies of an edge of $K_{n}$. The following definition generalizes this situation. 
\begin{madeflemma}[Symmetric extension]\label{lemmadefsymext}
Let $\L$ be a geometric lattice, $(\E, \imath)$ a modular extension of $\L$ and $F$ a modular coatom of $\L$. The subset 
$$\M \coloneqq \{F \in \L\cup_{[\zero, F]}\E \,\,  | \, \,\exists H \in \At(\L)\setminus [\zero, F] \, \textrm{ s.t. }\, (H, \imath(H))\leq F\}$$
is a modular cut of $\L\cup_{[\zero, F]}\E$. The associated single-element extension will be denoted by $\L\cup_F\E\cup_Fe$ and will be called the \textit{symmetric extension along $F$}. 
\end{madeflemma}
\begin{ex}
Let us go back to the geometric lattice $\L'$ described in \ref{exgeolatt}. This geometric lattice is isomorphic to the symmetric extension of the modular extension $\L = \L'_{|\{a,b,c,d\}} \hookrightarrow \E = \L$ along the modular flat $\{a\}$. We draw below the corresponding point-line configurations. 
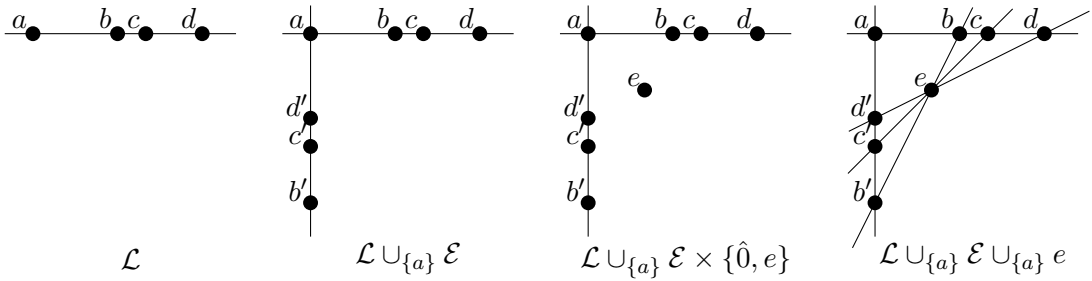
\begin{equation*}
    \begin{tikzpicture}[scale= 0.5, baseline=0pt]
    \node[circle, fill, black, inner xsep=0.7mm,inner ysep = 0.7mm]  (A) at (0,0) {} node[above left =-5pt of A, outer sep=0pt] {$a$};
    \node[circle, fill, black, inner xsep=0.7mm,inner ysep = 0.7mm]  (B) at (3,0) {};
    \node[circle, fill, black, inner xsep=0.7mm,inner ysep = 0.7mm]  (E) at (2.25,0) {};
    \node[circle, fill, black, inner xsep=0.7mm,inner ysep = 0.7mm]  (I) at (4.5,0) {};
    \draw[black] ($ (A)!-0.25!(B) $) -- ($ (A)!1.8!(B) $); 
    \node[above left =-5pt of E, outer sep=0pt] {$b$};
    \node[above left =-5pt of B, outer sep=0pt] {$c$};
    \node[above left =-5pt of I, outer sep=0pt] {$d$};
    \node[] at (2.6,-6) {$\L$};
    \end{tikzpicture}
    \hspace{15pt}
    \begin{tikzpicture}[scale= 0.5,baseline=0pt]
    \node[circle, fill, black, inner xsep=0.7mm,inner ysep = 0.7mm]  (A) at (0,0) {} node[above left =-5pt of A, outer sep=0pt] {$a$};
    \node[circle, fill, black, inner xsep=0.7mm,inner ysep = 0.7mm]  (B) at (3,0) {};
    \node[circle, fill, black, inner xsep=0.7mm,inner ysep = 0.7mm]  (C) at (0,-3) {};
    \node[circle, fill, black, inner xsep=0.7mm,inner ysep = 0.7mm]  (E) at (2.25,0) {};
    \node[circle, fill, black, inner xsep=0.7mm,inner ysep = 0.7mm]  (F) at (0,-2.25) {};
    \node[circle, fill, black, inner xsep=0.7mm,inner ysep = 0.7mm]  (I) at (4.5,0) {};
    \node[circle, fill, black, inner xsep=0.7mm,inner ysep = 0.7mm]  (J) at (0,-4.5) {};
    \node[above left =-5pt of E, outer sep=0pt] {$b$};
    \node[above left =-5pt of B, outer sep=0pt] {$c$};
    \node[above left =-5pt of I, outer sep=0pt] {$d$};
    \node[above left =-7pt of F, outer sep=0pt] {$d'$};
    \node[above left =-7pt of C, outer sep=0pt] {$c'$};
    \node[above left =-7pt of J, outer sep=0pt] {$b'$};
    \draw[black] ($ (A)!-0.25!(B) $) -- ($ (A)!1.8!(B) $); 
    \draw[black] ($ (A)!-0.25!(C) $) -- ($ (A)!1.8!(C) $);   
    \node[] at (2.6,-6) {$\L\cup_{\{a\}}\E $};
    \end{tikzpicture}
    \hspace{15pt}
    \begin{tikzpicture}[scale= 0.5,baseline=0pt]
    \node[circle, fill, black, inner xsep=0.7mm,inner ysep = 0.7mm]  (A) at (0,0) {} node[above left =-5pt of A, outer sep=0pt] {$a$};
    \node[circle, fill, black, inner xsep=0.7mm,inner ysep = 0.7mm]  (B) at (3,0) {};
    \node[circle, fill, black, inner xsep=0.7mm,inner ysep = 0.7mm]  (C) at (0,-3) {};
    \node[circle, fill, black, inner xsep=0.7mm,inner ysep = 0.7mm]  (D) at (1.5,-1.5) {};\node[circle, fill, black, inner xsep=0.7mm,inner ysep = 0.7mm]  (E) at (2.25,0) {};
    \node[circle, fill, black, inner xsep=0.7mm,inner ysep = 0.7mm]  (F) at (0,-2.25) {};
    \node[circle, fill, black, inner xsep=0.7mm,inner ysep = 0.7mm]  (I) at (4.5,0) {};
    \node[circle, fill, black, inner xsep=0.7mm,inner ysep = 0.7mm]  (J) at (0,-4.5) {};
    \node[above left =-5pt of E, outer sep=0pt] {$b$};
    \node[above left =-5pt of B, outer sep=0pt] {$c$};
    \node[above left =-5pt of I, outer sep=0pt] {$d$};
    \node[above left =-7pt of F, outer sep=0pt] {$d'$};
    \node[above left =-7pt of C, outer sep=0pt] {$c'$};
    \node[above left =-7pt of J, outer sep=0pt] {$b'$};
    \node[above left =-6pt of D, outer sep=0pt] {$e$};
    \draw[black] ($ (A)!-0.25!(B) $) -- ($ (A)!1.8!(B) $); 
    \draw[black] ($ (A)!-0.25!(C) $) -- ($ (A)!1.8!(C) $);
    
    \node[] at (2.6,-6) {$\L\cup_{\{a\}}\E\times \{\zero, e\}$};
    \end{tikzpicture}
    \hspace{15pt}
    \begin{tikzpicture}[scale= 0.5,baseline=0pt]
    \node[circle, fill, black, inner xsep=0.7mm,inner ysep = 0.7mm]  (A) at (0,0) {} node[above left =-5pt of A, outer sep=0pt] {$a$};
    \node[circle, fill, black, inner xsep=0.7mm,inner ysep = 0.7mm]  (B) at (3,0) {};
    \node[circle, fill, black, inner xsep=0.7mm,inner ysep = 0.7mm]  (C) at (0,-3) {};
    \node[circle, fill, black, inner xsep=0.7mm,inner ysep = 0.7mm]  (D) at (1.5,-1.5) {};\node[circle, fill, black, inner xsep=0.7mm,inner ysep = 0.7mm]  (E) at (2.25,0) {};
    \node[circle, fill, black, inner xsep=0.7mm,inner ysep = 0.7mm]  (F) at (0,-2.25) {};
    \node[circle, fill, black, inner xsep=0.7mm,inner ysep = 0.7mm]  (I) at (4.5,0) {};
    \node[circle, fill, black, inner xsep=0.7mm,inner ysep = 0.7mm]  (J) at (0,-4.5) {};
    \node[above left =-5pt of E, outer sep=0pt] {$b$};
    \node[above left =-5pt of B, outer sep=0pt] {$c$};
    \node[above left =-5pt of I, outer sep=0pt] {$d$};
    \node[above left =-7pt of F, outer sep=0pt] {$d'$};
    \node[above left =-7pt of C, outer sep=0pt] {$c'$};
    \node[above left =-7pt of J, outer sep=0pt] {$b'$};
    \node[above left =-6pt of D, outer sep=0pt] {$e$};
    \draw[black] ($ (A)!-0.25!(B) $) -- ($ (A)!1.8!(B) $); 
    \draw[black] ($ (A)!-0.25!(C) $) -- ($ (A)!1.8!(C) $);
    \draw[black] ($ (B)!-0.22!(C) $) -- ($ (B)!1.24!(C) $);
    \draw[black] ($ (E)!-0.465!(D) $) -- ($ (E)!3.8!(D) $);
    \draw[black] ($ (F)!-0.465!(D) $) -- ($ (F)!3.8!(D) $);
    \node[] at (2.6,-6) {$\L\cup_{\{a\}}\E\cup_{\{a\}} e $};
    \end{tikzpicture}
\end{equation*}
\end{ex}
\begin{proof}
To simplify notations let us assume that $\imath$ is an inclusion. Let $L_1 = (F_1, F'_1)$ and $L_2 = (F_2, F'_2)$ be elements of $\M \subset \L\cup_{[\zero, F]}\E$ such that $L_1 \wedge L_2 \notin \M$. We have to prove that $(L_1, L_2)$ is not a modular pair. We invite the reader to follow the argument on Figure \ref{figproofmodcut}. By assumption there exists $H_1 \neq H_2 \in \At(\L)\setminus [\zero, F]$ such that we have $H_1 \leq F_1, \, H_1 \leq F'_1,\, H_2 \leq F_2$ and $H_2 \leq F'_2$. Let us denote $F''_1 = F_1\wedge F$ and $F''_2 = F_2 \wedge F$. Let us prove the equality $L_1 \vee L_2 = (H_1,\zero) \vee (F''_1,F_1) \vee (F''_2,F_2)$. By modularity of $(F,\un)$ in $\L\cup_F\E$ (Lemma \ref{lemmapushout} item 2) we have $\rk(L_1) = 1 + \rk(F''_1, F'_1)$ which implies that we have $L_1 = (H_1,\zero)\vee (F''_1, F'_1)$. For any atom $H'_2$ below $F_2$ and not below $F$, by modularity of $F$ in $\L$ the element $(H_2\vee H'_2)\wedge F$ is an atom, which is below $F''_2$, and so also below $F'_2$. However by assumption $H_2$ is also below $F_2$, and so $H'_2 \leq ((H_2\vee H'_2)\wedge F) \vee H_2$ is also below $F'_2$. By modularity of $F$ in $\E$ the element $(H_1\vee H'_2)\wedge F$ is an atom of $F$, which is below $F'_1\vee F'_2$. Consequently the atom $H'_2 \leq H_1\vee ((H_1\vee H'_2)\wedge F)$ viewed in $\L\cup_F\E$ is below $(H_1,\zero)\vee (F''_1, F'_1) \vee (F''_2, F'_2)$ which concludes the promised equality. This implies that we have $\rk(L_1\vee L_2) = \rk ((F''_1,F_1) \vee (F''_2,F_2)) + 1$. 
On the other hand we have $\rk(L_1) = \rk(F'_1) + 1$ and $\rk(L_2) = \rk(F'_2) + 1$, which finally gives 
\begin{align*}
    \rk(L_1) + \rk(L_2) &= \rk(F'_1) + \rk(F'_2) + 2 \\ 
    &\geq \rk(F'_1\vee F'_2) + \rk(F'_1\wedge F'_2) + 2 \\
    &\geq \rk(L_1\vee L_2) + \rk(F'_1\wedge F'_2) + 1. 
\end{align*}
On the other hand we have $F'_1\wedge F'_2 = L_1\wedge L_2$ because as we have seen above, for $k=1,2$ if an atom belongs to $F_k$ then it belongs to $F'_k$. 
\end{proof}
\begin{figure}[H]
    \centering
    \begin{tikzpicture}
        \node[circle, fill, black, inner sep = 1.5pt] (A) at (3,3) {};
        \node[circle, fill, black, inner sep = 1.5pt] (B) at (3,-3) {};
        \node[circle, fill, black, inner sep = 1.5pt] (C) at (1.4,-0.2) {};
        \node[circle, fill, black, inner sep = 1.5pt] (D) at (4.9,-0.1) {};
        \node[circle, fill, black, inner sep = 1.5pt] (E) at (3.8,0.1) {};
        \draw (3,0) ellipse (3cm and 0.5cm);
        \draw[] (A) -- (0,0) ; 
        \draw[] (A) -- (6,0) ; 
        \draw[] (B) -- (6,0) ; 
        \draw[] (B) -- (0,0) ; 
        \node[] at (2.84, 2) {$\L$}; 
        \node[] at (2.8, -2) {$\L$}; 
        \node[] at (5.6, 2.7) {$\E$};
        \node[] at (0.84, 0) {$F$};
        \draw [] plot [smooth] coordinates {(0, 0) (0.1,1.5) (-0.1,3) (3,4.5) (6.5,3) (5.5,1.5) (6, 0)};
        \draw [red] plot [smooth] coordinates {(3, 3) (7, 0) (3,-3)};
        \node[red] at (7.2, 0) {$e$};
        \draw[green] (A) -- (C) node[midway, left] {$H_1$}; 
        \draw[green] (B) -- (C) node[midway, left] {$H_1$}; 
        \draw[blue] (B) -- (D) node[midway, right] {$H_2$}; 
        \draw[blue] (A) -- (D) node[midway, right] {$H_2$}; 
        \draw[blue] (B) -- (E) node[midway, left] {$H'_2$}; 
        \draw[blue, dashed] (D) -- (E);
        \draw[violet, dashed] (C) -- (E);
        \draw[blue] (A) -- (E) node[midway, left] {$H'_2$};
        
    \end{tikzpicture}
    \caption{Proof of Lemma \ref{lemmadefsymext}}
    \label{figproofmodcut}
\end{figure}
In Figure \ref{figproofmodcut} the edges in green belong to $L_1$, the edges in blue belong to $L_2$, the edge in red is the atom we are trying to add to $\L\cup_{F}\E$ and the edge in violet only belongs to $L_1\vee L_2$. The dashed edges represent atoms whose existence come from the modularity of $F$. 
\begin{lemma}\label{lemmamodinsym}
The embedding $\L\hookrightarrow \L\cup_F\E\cup_Fe$ is a modular extension. 
\end{lemma}
\begin{proof}
To simplify notations let us assume that $\imath$ is an inclusion. We have to prove that the element $(\un, F)$ is a modular flat of $\L\cup_F\E\cup_F e$. The flats $L$ of $\L\cup_F\E\cup_Fe$ are of three types: first it is possible that $e$ does not belong to $L$, in which case $L$ can be identified with a flat of $\L\cup_F\E$ which is not in $\M$. In this case we have 
\begin{align*}
    \rk((\un, F)\vee L) + \rk((\un, F)\wedge L) &= \rk_{\L\cup_F \E}((\un, F)\vee L) + \rk_{\L\cup_F \E}((\un, F)\wedge L) \\
    &= \rk_{\L\cup_F \E}((\un, F)) + \rk_{\L\cup_F \E}(L) \\
    &= \rk((\un, F)) + \rk(L),
\end{align*}
where the second equality comes from Lemma \ref{lemmapushout} item 2). Second, it is possible that $L$ contains $e$ and $L\setminus e$ is not a flat of $\L\cup_F\E\cup_F e$, which means that $L\setminus e$ can be identified with a flat of $\L\cup_F\E$ in $\M$. In this case we have
\begin{align*}
    \rk((\un, F)\vee L) + \rk((\un, F)\wedge L) &= \rk_{\L\cup_F \E}((\un, F)\vee (L\setminus e)) + \rk_{\L\cup_F \E}((\un, F)\wedge (L\setminus e))  \\
    &= \rk_{\L\cup_F \E}((\un, F)) + \rk_{\L\cup_F \E}(L\setminus e) \\
    &= \rk((\un, F)) + \rk(L).
\end{align*}
Finally, it is possible that $e$ belongs to $L$ and $L\setminus e$ is a flat. By construction this implies that $L$ does not contain any atom of the form $(H,\zero)$ or $(\zero, H)$ where $H$ is an atom of $\L$ which is not below $F$. In this case we have 
\begin{align*}
    \rk((\un, F)\vee L) + \rk((\un, F)\wedge L) &= \rk_{\L\cup_F \E}((\un, F)\vee (L\setminus e)) +1 + \rk_{\L\cup_F \E}((\un, F)\wedge (L\setminus e))  \\
    &= \rk_{\L\cup_F \E}((\un, F)) + \rk_{\L\cup_F \E}(L\setminus e) + 1 \\
    &= \rk((\un, F)) + \rk(L).
\end{align*}
\end{proof}
We are now ready to prove Theorem \ref{theoqiso}. The proof is an adaptation of the proof of \cite[Theorem 8.1]{lambrechts_formality_2012}.
\begin{proof}
Let $\L$ be a supersolvable geometric lattice with maximal chain of modular elements $\zero = G_0< \cdots < G_n= \un$. First notice that $\H(\I)$ is surjective, so we only need to prove that $\H^k(\MD(\L))$ and $\OS^k(\L)$ have the same dimension for all $k$. Recall from Definition/Proposition \ref{defgrading} that $\MD(\L)$ is graded by $\L$. This is also the case of $\OS(\L)$ (see \cite[Proposition 2.3]{Yuzvinsky_2001} for instance) and $\I_{\L}$ is compatible with those $\L$-gradings. By Definition/Proposition \ref{defgrading} we only have to prove the equality of dimension in grading $\un$, which we can assume to be different from $\zero$ by Remark \ref{rmkgradingzero}. 

It is a classical fact about Orlik--Solomon algebras that $OS(\L, \un)$ is concentrated in degree $\rk(\L)$ (see \cite[Proposition 2.3]{Yuzvinsky_2001} for instance). For all $1 \leq i\leq n$, denote by $J_i$ the set of atoms of $\L$ which are below $G_i$ and not below $G_{i-1}$. It is a classical fact about supersolvable geometric lattices that $\OS(\L, \un)$ has dimension $\#J_1\#J_2\cdots \#J_n$, where $\#J_i$ denotes the cardinality of $J_i$. More generally the Hilbert-Poincaré series of $\OS(\L)$ is $\prod_{i\leq n}(1+ \#J_it)$ (see \cite[Theorem 2]{Stanley1971} for instance). Let us prove by induction on $n$ that $H^{\bullet}(\MD(\L, \un))$ is concentrated in degree $n$ and has dimension $\#J_1\#J_2\cdots \#J_n$. Notice that the subcomplexes $H\bullet \MD(\L, G_{n-1})\subset \MD(\L, \un)$ for $H$ in $J_n$ are in direct sum. Let us denote $\MD^{\mathrm{top}}(\L) \coloneqq \bigoplus_{H \in J_n}H\bullet \MD(\L, G_{n-1})\subset \MD(\L, \un)$. By our induction hypothesis the cohomology of this complex is concentrated in degree $n$ and has dimension $\#J_1\#J_2\cdots \#J_n$. What remains is proving that the quotient complex $\MD(\L, \un)/ \MD^{\mathrm{top}}(\L)$ is acyclic. As a graded vector space the quotient complex $\Qu \coloneqq \MD(\L, \un)/ \MD^{\mathrm{top}}(\L)$ decomposes as a direct sum $\Qu = \Qu_{1}\oplus \Qu_{2}$ where $\Qu_1$ is spanned by modular diagrams $(\E, \imath, J)$ satisfying the condition that there exists a modular coatom $F$ of $\E$ such that $F\wedge F_{\imath} = G_{n-1}$ and such that $J\cap F^c$ is a singleton not contained in $\At(\L)$, and $\Qu_2$ is spanned by all the remaining modular diagrams. We have the following lemma. 
\begin{lemma}\label{lemmasymext}
Let $(\E, \imath, J)$ be a modular diagram in $\Qu_1$, with modular coatom $F$ and unique atom $e\in J\cap F^c$. The restriction $\E_{|\At(\L)\cup F \cup e}$ is isomorphic to a symmetric extension $\L\cup_{G_{n-1}}[\zero,F]\cup_{G_{n-1}} e$. 
\end{lemma}
\begin{proof}[Proof of Theorem \ref{theoqiso}]
Recall from Lemma \ref{lemmaisodiamond} that taking the join with $e$ gives an isomorphism $\vee_{e,F}:[\zero,F]\xrightarrow{\sim}[e,\un]$ and also gives an isomorphism $\vee_{e,F_{\imath}}:[\zero, F_{\imath}] \xrightarrow{\sim} [e, e\vee F_{\imath}]$. As a consequence the composition $\vee_{e,F}^{-1}\circ\vee_{e,F_{\imath}}\circ \imath$ defines an embedding of $\L$ in $[\zero,F]$ and we can consider the symmetric extension $\L\cup_{G_{n-1}}[\zero,F]\cup_{G_{n-1}} e$ along this embedding. By Lemma \ref{lemmadetectionpushout} the restriction $\E_{|\At(\L)\cup F}$ is isomorphic to the pushout $\L\cup_{G_{n-1}}[\zero,F]$. We have to show that a set $J$ of atoms in $\At(\L)\cup F$ generates $e$ if and only if $J$ contains the two copies $\imath(H)$ and $\vee_{e,F}^{-1}(e\vee\imath(H))$ of some atom $H$ of $\L$. By definition, for any atom $H$ of $\L$ the atom $\vee_{e,F}^{-1}(e\vee\imath(H))$ is below the rank $2$ flat $e\vee\imath(H)$, which means that we have $e \in \imath(H) \vee \vee_{e,F}^{-1}(e\vee\imath(H))$. For the other direction, let $F$ be a flat of $\E_{|\At(\L)\cup F}$ which is not a flat of $\E_{|\At(\L)\cup F\cup e }$. This implies that $F$ contains an atom which is not in $F$, say $\imath(H)$. By the previous considerations the elements $e$ and $\imath(H)$ generate $\vee_{e,F}^{-1}(e\vee\imath(H))$. This implies that the flat $F\cup e $ must contain $\vee_{e,F}^{-1}(e\vee\imath(H))$, and so does $F$.  
\end{proof}
Let us introduce a filtration on $\Qu$ indexed by the cardinality of $J$ on $\Qu_1$ and the cardinality of $J$ minus 1 on $\Qu_2$. We aim to prove that the associated graded of this filtration is acyclic. We have the following lemma.  
\begin{lemma}
The differential of the associated graded of $\Qu$ is zero on $\Qu_2$ and sends a modular diagram $\Gamma \in \Qu_1$ with special element $e$ to $\Gamma/e \in \Qu_2$.
\end{lemma}
\begin{proof}
The first statement is implied by the fact that the differential on $\MD(\L)$ strictly lowers the cardinality of $J$. Let $\Gamma = (\E, \imath, J)$ be a modular diagram in $\Qu_1$, with special element $e$. If $H$ is an atom in $J\setminus e$, notice that $F$ is a modular coatom of $[H, \un]$ and $H\vee e$ is the only atom of $H\vee J$ outside of $F$ in $[H, \un]$. In other words $\Gamma/H$  belongs to $\Qu_1$ and once again $H\vee J$ has cardinality one less than $J$, so $\Gamma/H$ must be zero in the associated graded. On the other hand, let us show that $\Gamma/e$ belongs to $\Qu_2$. By Lemma \ref{lemmasymext} we can assume that $\E$ is $\L\cup_{G_{n-1}}[\zero, F]\cup_{G_{n-1}} e$. If $[e,\un]$ had a modular coatom $F'$ with a unique element $e\vee H \in e\vee J$ outside of $F'$. By Lemma \ref{lemmaisodiamond} the flat $\vee_{e,F}^{-1}(F')$ is modular in $[\zero, F]$ and so by Lemma \ref{lemmamodularpushout} the flat $F_{\imath}\vee \vee_{e,F}^{-1}(F')$ is modular in $\L\cup_{G_{n-1}}[\zero, F]$. By the same arguments as in the proof of Lemma \ref{lemmamodinsym}, the flat $F_{\imath}\vee \vee_{e,F}^{-1}(F')$ is also modular in $\L\cup_{G_{n-1}}[\zero, F]\times _{G_{n-1}}e = \E$. However, there are only two atoms of $J$ outside of $F_{\imath}\vee \vee_{e,F}^{-1}(F')$, namely $e$ and $H$, and so by relation \ref{rel4} in Definition \ref{defMD} the modular diagram $\Gamma$ is zero in $\Qu$. 
\end{proof}
By Lemma \ref{lemmasymext} the map $\Gamma \rightarrow \Gamma/e$ is an isomorphism with inverse given by sending a diagram $(\E, \imath, J) \in \Qu_2$ to $(\L\cup_{G_{n-1}}\E\cup_{G_{n-1}}e, \imath, J\cup e)\in \Qu_1$. By a standard spectral sequence argument this implies that $\Qu$ is acyclic, which concludes the proof of Theorem \ref{theoqiso}. 
\end{proof}
\section{Additional structures on $\MD(\L)$}\label{secstruct}
In this section we introduce additional structures on $\MD(\L)$, which generalize the additional structures on the cdga of admissible graphs described in Section \ref{secgra}. First notice that for any geometric lattices $\L$ and $\L'$, an isomorphism $\varphi:\L\xrightarrow{\sim} \L'$ induces an isomorphism $\MD(\varphi): \MD(\L')\xrightarrow{\sim}\MD(\L)$ defined by  
\begin{equation*}
    \MD(\varphi)(\E, \imath, J) \coloneqq (\E, \imath\circ \varphi, J).
\end{equation*}
One can see that this is an isomorphism of cdgas, and the map $\varphi \rightarrow \MD(\varphi)$ is compatible with composition on both sides. For any partition lattice $\Pi_n$ with $n\geq 3$, the automorphism group of $\Pi_n$ is $\Sym_n$ and we recover the classical symmetric action. For $n=2$ however, the automorphism group of $\Pi_2 = \{\zero < \un \}$ is the trivial group instead of $\Sym_2$ and so we do not recover the $\Sym_2$ action. This is not surprising because $\Sym_2$ acts trivially on $\H^{\bullet}(\S^1)$ anyway. Note that we have a similar equivariant structure for Orlik--Solomon algebras, with isomorphisms $\OS(\varphi)$ defined by $\OS(\varphi)(e_J) \coloneqq e_{\varphi^{-1}(J)}.$ The morphisms $\I_{\L}$ are compatible with those equivariant structures in an obvious way.\\

As we saw in Section \ref{secgra}, the cdga of admissible graphs admits a (co)operadic structure, which models the inclusions of the strata of the Fulton--Macpherson compactification of the braid arrangements. In \cite{BC2024} and \cite{BC2025} it was explained by the author how to construct a type of operad-like structure mimicking inclusions of strata of compactifications of hyperplane arrangement complements, in a purely combinatorial setting. Let us recall the main definition here. 
\begin{madef}[$\GL$-cooperad]
A $\GL$-\textit{cooperad} $\C$ in a monoidal category $(\mathfrak{C},\otimes)$ is a collection of objects $\C(\L) \in \mathfrak{C}$ indexed by non-trivial geometric lattices, together with a collection of morphisms $\Delta_{F,\L}:  \C(\L)\rightarrow \C([\zero, F])\otimes \C([F, \un]) $ indexed by any choice of a proper flat $F$ in a geometric lattice $\L$, satisfying the condition 
\begin{equation}\label{eqdefcoop} 
(\Delta_{F_1, [\zero, F_2]} \otimes \Id_{\C([F_2, \un])})\circ \Delta_{F_2, \L}= (\Id_{\C([\zero, F_1])} \otimes \Delta_{F_2, [F_1, \un]} ) \circ \Delta_{F_1, \L}
\end{equation}
for all geometric lattice $\L$ and all proper flats $F_1 < F_2$ of $\L$. 
\end{madef}
We will simply write $\Delta_F$ if the geometric lattice $\L$ can be deduced from the context. A $\GL$-cooperad restricted to partition lattices is a shuffle cooperad with levels (see \cite{DK_2010} and \cite{Fresse_2003}). One can get rid of the levels by introducing the notion of building sets and use partition lattices with minimal building sets instead of maximal building sets (see \cite{BC2024} for more details). 
\begin{ex}[{\cite[Example 2.2]{BC2025}}]\label{exOS}
The collection of Orlik--Solomon algebras has a $\GL$-cooperadic structure in the monoidal category of graded commutative algebras, with cooperadic coproducts $\Delta_{F}$ defined by sending a generator $e_H$ to $e_H\otimes 1$ if $H$ is below $F$ and to $1\otimes e_{F\vee H}$ otherwise. This cooperadic structure is compatible with the equivariant structure on $\OS$ in the sense that for any isomorphism $\L_1\xrightarrow{\sim} \L_2$ and any flat $F \in \L$ we have the equality of morphisms $\Delta_{F, \L}\circ \OS(\varphi) = (\OS(\varphi_{|[\zero, F]}) \otimes \OS(\varphi_{|[F,\un]}))\circ \Delta_{\varphi(F), \L'}.$ 
\end{ex}
We also have a natural notion of morphism between $\GL$-cooperads.
\begin{madef}[Morphism of $\GL$-cooperads]
Let $(\C_1, \Delta_1)$ and $(\C_2, \Delta_2)$ be two $\GL$-cooperads in some monoidal category $(\mathfrak{C}, \otimes)$. A \textit{morphism of $\GL$-cooperads} $f:\C_1\rightarrow \C_2$ is a collection of morphims $f_{\L}: \C_1(\L)\rightarrow \C_2(\L)$ in $\mathfrak{C}$ indexed by nontrivial geometric lattices, satisfying the equality $(\Delta_2)_{F, \L} \circ f_{\L} = (f_{[\zero, F]}\otimes f_{[F,\un]})\circ (\Delta_1)_{F,\L}$ for all proper flat $F \in \L$. 
\end{madef}
\begin{prop}\label{propcoop}
The collection $\MD$ admits a $\GL$-cooperadic structure in cdgas, and $\I$ is a morphism of $\GL$-cooperads between $\MD$ and $\OS$. 
\end{prop}
\begin{proof}
Let $\Gamma = (\E, \imath, J)$ be a modular diagram of some geometric lattice $\L$. For any proper flat $F \in \L$ and any flat $F'\in \E$ such that $F_{\imath} \wedge F' = F $ we put 
\begin{align*}
    \Gamma_{F'} &\coloneqq ([\zero, F'], \imath_{|[\zero, F]}, J \cap F'), \\ 
    \Gamma^{F'} &\coloneqq ([F',\un], \vee_{F, F_{\imath}} \circ \imath_{|[F,\un]}, F'\vee (J\cap F'^c)),
\end{align*}
which are modular diagrams of $[\zero, F]$ and $[F,\un]$ respectively, by Lemmas \ref{lemmamodheredity} and \ref{lemmaisodiamond}. We also define $\epsilon(J,F')$ to be the signature of the permutation of $J$ keeping the relative order in $J\cap F'$ and $J\cap F'^c$, and sending every element of $J\cap F'$ before every element of $J\cap F'^c$. We then define the cooperadic coproduct $\Delta_{F, \L}: \MD(\L) \rightarrow \MD([\zero, F])\otimes \MD([F, \un])$ by 
\begin{equation*}
    \Delta_{F,\L}(\Gamma) \coloneqq \sum_{\substack{F' \in \E \textrm{ s.t. }\\ \imath^{-1}(F'\wedge F_{\imath}) = F }}\epsilon(J,F')\Gamma_{F'}\otimes \Gamma^{F'}
\end{equation*}
on modular diagrams of $\L$. Let us quickly explain why $\Delta_{F, \L}$ is compatible with relations \ref{rel1} to \ref{rel5} in Definition \ref{defMD}. For relation \ref{rel1} any exchange of two elements of $J\cap F'$ or two elements of $J\cap F'^c$ will induce a sign change in $\Gamma_{F'}$ and $\Gamma^{F'}$ respectively, and an exchange of two elements one in $J\cap F'$ and the other in $J\cap F'^c$ will induce a sign change in $\epsilon(J,F')$. For relation \ref{rel2}, if we have an embedding $\imath':\E \hookrightarrow \E'$ one can see that the flats of $\E'$ whose intersection with $F_{\imath'\circ \imath}$ is $\imath'\circ \imath (F)$ are exactly the flats of the form $\imath'(F')$ with $F'\in \E$ having intersection $\imath(F)$ with $F_{\imath}$. This means that the terms in $\Delta_{F,\L}((\E, \imath, J))$ and $\Delta_{F,\L}((\E, \imath'\circ\imath, \imath'(J))$ respectively are in bijection and one can check that those terms coincide. For relation \ref{rel3}, if we have $F_{\imath}\vee \bigvee J< \un$ then $\Gamma_{F'}$ is zero if $F'$ is not below $F_{\imath}\vee \bigvee J$ and $\Gamma^{F'}$ is zero if $F'$ is below $F_{\imath}\vee \bigvee J$. For relation \ref{rel4}, if $\E$ is a non-trivial product $\E_1\times \E_2$ with $\E_1$ containing $F_{\imath}$, then one can see that for any $F'$, at least one of $[\zero, F']$ and $[F', \un]$ will be a non-trivial product. For relation \ref{rel5} if there exists a modular coatom not containing only two elements $H_1, H_2 \in J_{\Gamma}$ then for all flat $F'$ either $\Gamma_{F'}$ or $\Gamma^{F'}$ is zero depending on if $F'$ contains $H_1, H_2$ or not. \\

We now prove that $\Delta_{F}$ is a morphism of algebras for all proper flat $F$. For all modular diagrams $\Gamma_1 = (\E_1, \imath_1, J_1)$ and $\Gamma_2 = (\E_2, \imath_2, J_2)$ of $\L$, we have 
\begin{align*}
    \Delta_{F}(\Gamma_1\bullet \Gamma_2) &= \sum_{\substack{F' \in \E_1\cup_{\L} \E_2 \\ \imath_{12}^{-1}(F' \wedge F_{\imath_{12}}) = F}}\epsilon(J_1\cup J_2,F')(\Gamma_1\bullet\Gamma_2)_{F'}\otimes (\Gamma_1\bullet\Gamma_2)^{F'} \\ 
    &= \sum_{\substack{F'_1 \in \E_1 \\ \imath_{1}^{-1}(F'_1 \wedge F_{\imath_{1}}) = F}}\sum_{\substack{F'_2 \in \E_2  \\\imath_{2}^{-1}(F'_2 \wedge F_{\imath_{2}}) = F}} \epsilon(J_1\cup J_2,(F'_1, F'_2))(\Gamma_1\bullet\Gamma_2)_{(F'_1, F'_2)}\otimes (\Gamma_1\bullet\Gamma_2)^{(F'_1, F'_2)} \\ 
    &= \sum_{\substack{F'_1 \in \E_1 \\ \imath_{1}^{-1}(F'_1 \wedge F_{\imath_{1}}) = F}}\sum_{\substack{F'_2 \in \E_2  \\\imath_{2}^{-1}(F'_2 \wedge F_{\imath_{2}}) = F}} \epsilon(J_1,F'_1)\epsilon(J_2, F'_2)((\Gamma_1)_{F'_1}\otimes (\Gamma_1)^{F'_1})\bullet ((\Gamma_2)_{F'_2}\otimes (\Gamma_2)^{F'_2}) \\ 
    &= \Delta_{F}(\Gamma_1)\bullet \Delta_{F}(\Gamma_2).
\end{align*} 

We now prove that the morphisms $\Delta_{F}$ satisfy equation \eqref{eqdefcoop}. Let $F_1<F_2$ be two proper flats of $\L$ and let $\Gamma$ be a modular diagram of $\L$. We have 
\begin{align*}
(\Delta_{F_1, [\zero, F_2]} \otimes \Id)(\Delta_{F_2, \L}(\Gamma)) &= (\Delta_{F_1, [\zero, F_2]} \otimes \Id)\left( \sum_{\substack{F' \in \E \textrm{ s.t. }\\ F'\wedge F_{\imath} = F_2 }}\epsilon(J,F')\Gamma_{F'}\otimes \Gamma^{F'}\right) \\
& = \sum_{\substack{F' \in \E \textrm{ s.t. }\\ F'\wedge F_{\imath} = F_2 }} \sum_{\substack{F'' \in [\zero, F'] \textrm{ s.t. }\\F'' \wedge F_{\imath} = F_1 }} \epsilon(J,F')\epsilon(J\cap F',F'') (\Gamma_{F'})_{F''}\otimes (\Gamma_{F'})^{F''} \otimes \Gamma^{F'}\\
&=\sum_{\substack{F'' \in \E \textrm{ s.t. }\\ F''\wedge F_{\imath} = F_1 }} \sum_{\substack{F' \in [F'', \un] \textrm{ s.t. }\\F' \wedge F_{\imath} = F_2 }} \epsilon(J,F')\epsilon(J\cap F',F'') \Gamma_{F''}\otimes (\Gamma^{F''})_{F'}\otimes (\Gamma^{F''})^{F'} \\ 
&= (\Id \otimes \Delta_{F_2, [F_1, \un]})\left( \sum_{\substack{F'' \in \E \textrm{ s.t. }\\ F''\wedge F_{\imath} = F_1 }}\epsilon(J,F'')\Gamma_{F''}\otimes \Gamma^{F''} \right) \\ 
&= (\Id \otimes \Delta_{F_2, [F_1, \un]})(\Delta_{F_1, \L}(\Gamma)).
\end{align*}

We now prove that $\Delta_{F}$ is a morphism of dg complexes for all proper flat $F$. This is more subtle than the previous results because this statement is not true before quotienting by relations \ref{rel1} to \ref{rel5}. Let $\Gamma = (\E, \imath, J)$ be a modular diagram of $\L$, and $F$ a proper flat of $\L$. We say that a flat $F'\in \E$ is $(F,\Gamma)$-admissible if $\imath^{-1}(F'\wedge F_{\imath})= F$, no element of $J \cap F'$ is a bridge in $\Gamma_{F'}$ and $F'\vee \cdot $ is injective on $J\cap F'^c$. Clearly we have 
\begin{align*}
    \Delta_{F}(\Gamma) = \sum_{F' \, (F,\Gamma)-\textrm{adm}}\epsilon(J,F')\Gamma_{F'}\otimes \Gamma^{F'}.
\end{align*}
If $F'$ is $(F,\Gamma)$-admissible we say that an element $H \in J$ is $F'$-contractible if either $H\in F'$ and $H$ is contractible in $\Gamma_{F'}$, or $H\notin F'$ and $F'\vee H$ is contractible in $\Gamma^{F'}$. We have the equalities
\begin{align*}
    \d(\Delta_{F}(\Gamma)) &= \d\left(\sum_{F' \, (F, \Gamma)-\textrm{adm}}\epsilon(J,F')\Gamma_{F'}\otimes \Gamma^{F'}\right) \\ 
    &= \sum_{F' \, (F, \Gamma)-\textrm{adm}}\epsilon(J,F')(\d(\Gamma_{F'})\otimes \Gamma^{F'} + (-1)^{\deg(\Gamma_{F'})}\Gamma_{F'}\otimes \d(\Gamma^{F'})) \\ 
    &= \sum_{F' \, (F, \Gamma)-\textrm{adm}} \epsilon(J,F')\left(\sum_{\substack{H_i \in \Sc\cap F'\\}}(-1)^{i+1} \Gamma_{F'}/H_i\otimes \Gamma^{F'} \right.\\ 
    & \hspace{80pt}\left.+(-1)^{\deg{\Gamma_{F'}}}\sum_{\substack{H_i\in \Sc \cap F'^c \\ H_i \nleq F'\vee F_{\imath}}} (-1)^{i+1}\Gamma_{F'}\otimes \Gamma^{F'}/(F'\vee H_i)\right) \\ 
    &= \sum_{F' \, (F, \Gamma)-\textrm{adm}} \epsilon(J,F')\left(\sum_{\substack{H_i \in \Sc\cap F'\\}}(-1)^{i+1} (\Gamma/H_i)_{F'}\otimes (\Gamma/H_i)^{F'} \right.\\ 
    & \hspace{80pt}\left.+(-1)^{\deg{\Gamma_{F'}}}\sum_{\substack{H_i\in \Sc \cap F'^c \\ H_i \nleq F'\vee F_{\imath}}} (-1)^{i+1}(\Gamma/H_i)_{F'\vee H_i}\otimes (\Gamma/H_i)^{F'\vee H_i}\right) \\ 
    &= \sum_{H_i \in \Sc} (-1)^{i+1}\sum_{F'' \, (F,\Gamma/H_i)-\textrm{adm.}}\epsilon(J_{\Gamma/H_i}, F'')(\Gamma/H_i)_{F''}\otimes (\Gamma/H_i)^{F''} \\
    &= \Delta_{F} (\d(\Gamma)). 
\end{align*}

Finally, one can easily check that $\I$ is compatible with the cooperadic structures on both side because $\I(\E, \imath, J)$ is zero whenever $J$ is not contained in $\At(\L)$, and when $J$ is contained in $\At(\L)$ the cooperadic coproducts coincide. 
\end{proof}
Note that the cooperadic structure on $\MD$ is compatible with the equivariant structure on $\MD$, in the same sense as described in Example \ref{exOS} for Orlik--Solomon algebras. \\

In principle Proposition \ref{propcoop} together with Theorem \ref{theoqiso} should lead to a proof of the formality of the $\GL$-operad of Fulton--Macpherson compactifications restricted to supersolvable complex hyperplane arrangements, in analogy with the braid case, but it is not clear to the author how one would use this result. 

\section{Koszulness of Orlik--Solomon algebras}\label{seckos}
A graded commutative algebra is called $1$-generated if it is generated by elements of grading $1$. It is quadratic if relations between degree $1$ generators are generated by elements of grading $2$. One can keep going and ask if relations between relations are generated in grading $3$ and so on. If this is true at all levels we say that our commutative graded algebra is Koszul. There are several ways to formalize this, one of them is via the so-called Chevalley--Eilenberg complex. Define the Lie Koszul dual of a quadratic graded commutative algebra $A = \Lambda^{\bullet}V/R$ to be the graded Lie algebra $A^{!} \coloneqq \Lie[V[-1]]/R^{\bot}$. Set also $A^{\text{!`}} \coloneq (A^{!})^{\vee}.$ The Lie Koszul dual of (the quadratic envelope of) an Orlik--Solomon algebra $\OS(\L)$ is the holonomy Lie algebra $\mathfrak{h}(\L)$ studied in \cite{guo2023}. It has the following presentation.
\begin{equation*}
    \mathfrak{h}(\L) = \frac{\Lie[t_{H}\, | \, H \in \At(\L), \deg(t_H) = 0]}{\langle [t_H, \sum_{H' \in F}t_{H'}], H \leq  F, \rk(F) = 2\rangle}.
\end{equation*}

\begin{madef}[Koszulness]
A quadratic graded commutative algebra $A$ is \textit{Koszul} if there exists a zig-zag of quasi-isomorphisms between $A$ and the Chevalley-Eilenberg complex $C^{\bullet}_{CE}(A^{\text{!`}}).$
\end{madef}
We have the following theorem. 
\begin{theo}[{\cite[Corollary~6.21]{Yuzvinsky_2001}}]
If a geometric lattice $\L$ is supersolvable, then $\OS(\L)$ is Koszul. 
\end{theo}
This result can be proved by exhibiting a quadratic Gröbner basis of $\OS(\L)$. We give here a new proof using the cdga of modular diagrams, which follows a strategy of proof suggested in \cite{SW_2009} in the braid case. 
\begin{proof}
By Theorem \ref{theoqiso}, if $\L$ is supersolvable we have a quasi-isomorphism $\I_{\L}: \MD(\L) \xrightarrow{\sim} \OS(\L)$. By Proposition \ref{propfree} the cdga $\MD(\L)$ is quasi-free with space of generators $\IMD(\L)$ with basis consisting of equivalence classes of irreducible modular diagrams. The differential on $\MD(\L)$ gives a structure of $L_{\infty}$-coalgebra on $\IMD(\L)$ and by definition we have $C^{\bullet}_{CE}(\IMD(\L)) = \MD(\L).$ Since the Chevalley--Eilenberg complex is functorial and preserves quasi-isomorphisms of $L_{\infty}$-coalgebras, what remains is proving that we have a quasi-isomorphism of $L_{\infty}$-coalgebras $\h(\L)^{\vee} \xrightarrow{\sim}\IMD(\L)$. The proof is a direct adaptation of that given in \cite[Appendix B]{SW_2009}. Let us repeat the argument here. Recall that if we have two sub Lie algebras $L_1, L_2$ of a Lie algebra $L$, we say that $L$ is a semi-direct product of $L_1$ and $L_2$, denoted $L = L_1 \ltimes L_2 $, if we have $L = L_1\oplus L_2$ as vector spaces, and if in addition we have $[L_1, L_2] \subset L_2.$ There is an analogous notion for $L_{\infty}$-algebras. Let $\zero = G_0 < G_1 < \cdots < G_{n-1} < G_n = \un$ be a maximal chain of modular elements in $\L.$ Recall from Definition \ref{defgrading} that we have an $\L$-grading on $\MD(\L)$, defined by $\gr (\E, \imath, J) = (\bigvee J) \wedge F_{\imath}$. Denote by $\F^{\vee}\subset \IMD(\L)^{\vee}$ the sub-$L_{\infty}$-algebra of irreducible modular diagrams with $\L$-grading not below $G_{n-1}.$ We have a semi-direct product $\IMD(\L)^{\vee} = \IMD([\zero, G_{n-1}])^{\vee} \ltimes \F^{\vee}$. On the other hand by \cite[Theorem 4.6]{guo2023} we have a semi-direct product $\h(\L) = \h([\zero,G_{n-1}])\ltimes \Lie[t_{H} \, | \, H \in J_n]$, where $J_n$ denotes the set of atoms of $\L$ which are not below $G_{n-1}.$ By induction we only need to prove that we have a quasi-isomorphism of $L_{\infty}$-coalgebras $\Phi: \F \xrightarrow{\sim} \Lie[t_{H} \, | \, H \in J_n]^{\vee}.$ Let us define $\Phi$ by induction. As a base case, for all $H\in J_n$ we set $\Phi(t_H^{\vee}) = \{H\}$. In general for some element $\alpha \in \Lie[J_n]$ with cobracket $\Delta(\alpha) = \sum_{i}\alpha^{(1)}_i\otimes \alpha^{(2)}$, with $\alpha^{(1)}_i = (\E^{(1)}_i, \imath_i^{(1)}, J_i^{(1)})$ and $\alpha^{(2)}_i = (\E^{(2)}_i, \imath_i^{(2)}, J_i^{(2)})$ we set $\Phi(\alpha) = \sum_i (\L\cup_{G_{n-1}}(\E^{(1)}_i\cup \E^{(2)}_i)\cup_{G_{n-1}} e, \imath, J^{(1)}_i\cup J^{(2)}_i).$ This is a morphism of $L_{\infty}$-coalgebras, which is injective in cohomology. This means we only have to prove that $\H(\F)$ has the right dimension. As in the proof of Theorem \ref{theoqiso} let us denote by $\mathcal{IQ}_1$ the subspace of $\F$ spanned by irreducible modular diagrams $\Gamma$ such that there exists a modular flat $F$ above $G_{n-1}$ satisfying the condition that $J_{\Gamma}$ contains a unique atom not below $F$. We also denote by $\mathcal{IQ}_2$ the span of all the other irreducible modular diagrams. The complex $\F$ is the total complex of the bi-complex 
\begin{equation*}
    \begin{tikzpicture}
        \node[] (A) at (0,0) {$\mathcal{IQ}^{i}_2$} ; 
        \node[] (B) at (-2,0) {$\mathcal{IQ}^{i-1}_2$} ; 
        \node[] (C) at (2,0) {$\mathcal{IQ}^{i+1}_2$} ; 
        \node[] (D) at (-2,-2) {$\mathcal{IQ}^{i-1}_1$} ; 
        \node[] (E) at (0,-2) {$\mathcal{IQ}^{i}_1$} ; 
        \node[] (F) at (2,-2) {$\mathcal{IQ}^{i+1}_1$} ; 
        \node[] (G) at (-4,0) {$\cdots$} ; 
        \node[] (H) at (4,-0) {$\cdots$} ; 
        \node[] (I) at (-4,-2) {$\cdots$} ; 
        \node[] (J) at (4,-2) {$\cdots$} ; 
        \draw[->] (B) -- (A); 
        \draw[->] (A) -- (C); 
        \draw[->] (D) -- (E); 
        \draw[->] (E) -- (F); 
        \draw[->] (D) -- (B) node[midway, left]{$/e$}; 
        \draw[->] (E) -- (A)node[midway, left]{$/e$}; 
        \draw[->] (F) -- (C) node[midway, left]{$/e$}; 
        \draw[->] (G) -- (B);
        \draw[->] (C) -- (H);
        \draw[->] (I) -- (D);
        \draw[->] (F) -- (J);

    \end{tikzpicture}
\end{equation*}
Let us consider the spectral sequence associated to considering first the vertical differentials and then the horizontal ones. As explained in the proof of Theorem \ref{theoqiso} the vertical arrows are surjective, with kernel the dg-complex $\mathcal{IQ}_1^{disc}$ spanned by diagrams in $\mathcal{IQ}_1$ whose contraction at $e$ is not irreducible. We can identify $\mathcal{IQ}_1^{disc}$ with $C^{\geq 2}_{CE}\F$, by sending a diagram of $\mathcal{IQ}_1^{disc}$ to the tensor product of the components after contracting $e$. Consider now the filtration given by the number of those irreducible components. The differential of the associated graded is simply the differential of $\F$. By induction on the cardinality of $J_{\Gamma}$ the cohomology of that differential is $C^{\geq 2}_{CE}\Lie[t_H\, | \, H\in J_n]$ and the second differential is the $CE$-differential. The cohomology of $C^{\geq 2}_{CE}\Lie[t_H\, | \, H\in J_n]$ with the $CE$-differential is $\Lie[t_H\, | \, H\in J_n]$ and we are done.
\end{proof}
To finish let us mention the following classical conjecture. 
\begin{conj}[{\cite[Section 6]{Yuzvinsky_2001}}]
The Orlik--Solomon algebra $\OS(\L)$ of a geometric lattice $\L$ is Koszul if and only if $\L$ is supersolvable.
\end{conj}
We hope this article could shed new light on that problem. We propose the following refinement. 
\begin{conj}
Let $\L$ be a geometric lattice. The following statements are equivalent: 
\begin{enumerate}
    \item $\L$ is supersolvable.
    \item $\OS(\L)$ is Koszul.
    \item $\I_{\L}$ is a quasi-isomorphism. 
\end{enumerate}
\end{conj}
\bibliographystyle{alpha}
\bibliography{sample}

\end{document}